\documentclass[reqno,11pt]{amsart}
\usepackage{relsize}
\usepackage{scalerel}
\usepackage[margin=1in]{geometry}
\usepackage{stackengine,wasysym}
\usepackage{todonotes}
\usepackage{xcolor}
\usepackage{mathrsfs}
\usepackage{dsfont}
\usepackage{mathtools}
\mathtoolsset{showonlyrefs}
\usepackage[hyperfootnotes=false,colorlinks=true,linkcolor = blue, urlcolor  = blue, citecolor = blue]{hyperref}
\usepackage[sort,nocompress]{cite}
\usepackage{float}
\usepackage{amsmath, amsthm, amssymb}
\usepackage{times}
\usepackage{color}
\usepackage{comment}
\usepackage[toc,page]{appendix}
\usepackage{verbatim}
\usepackage{enumerate}

\usepackage{tikz}

\newcommand{\pa}{\partial}
\newcommand{\la}{\label}
\newcommand{\fr}{\frac}
\newcommand{\na}{\nabla}
\newcommand{\be}{\begin{equation}}
\newcommand{\ee}{\end{equation}}
\newcommand{\bes}{\begin{equation*}}
\newcommand{\ees}{\end{equation*}}
\newcommand{\ba}{\begin{array}{l}}
\newcommand{\ea}{\end{array}}

\newcommand{\beg}{\begin}

\newcommand{\veps}{\varepsilon}

\newcommand{\abs}[1]{\left\lvert#1\right\rvert}

\newcommand{\C}{\mathbb C}
\newcommand{\R}{\mathbb R}
\newcommand{\Q}{\mathbb Q}

\def\RR{{\mathbb R}}
\def\TT{{\mathbb T}}

\newtheorem{Thm}{Theorem}[section]%  meant for continuous numbers
\newtheorem{prop}[Thm]{Proposition}

\numberwithin{equation}{section}

\title{ERROR ESTIMATES OF PHYSICS-INFORMED NEURAL NETWORKS FOR APPROXIMATING  BOLTZMANN EQUATIONS}

\author[E. Abdo]{Elie Abdo}
\address[E. Abdo]
{	Department of Mathematics \\
	University of California  \\
	Santa Barbara, CA 93106-3080, USA.} \email{elieabdo@ucsb.edu}

\author[L. Chai]{Lihui Chai}
\address[L. Chai]
{	School of Mathematics \\
	Sun Yat-sen University \\
	Guangzhou, 510275, China.} \email{chailihui@mail.sysu.edu.cn}

\author[R. Hu]{Ruimeng Hu}
\address[R. Hu]
{Department of Mathematics \\
	Department of Statistics and Applied Probability \\
	University of California  \\
	Santa Barbara, CA 93106-3080, USA.} \email{rhu@ucsb.edu}

\author[X. Yang]{Xu Yang}
\address[X. Yang]
{Department of Mathematics \\
	University of California  \\
	Santa Barbara, CA 93106-3080, USA.} \email{xuyang@math.ucsb.edu}

\begin{document}
\begin{abstract}
		Motivated by the recent successful application of physics-informed neural networks (PINNs) to solve Boltzmann-type equations [S. Jin, Z. Ma, and K. Wu, J. Sci. Comput., 94 (2023), pp. 57], we provide a rigorous error analysis for PINNs in approximating the solution of the Boltzmann equation near a global Maxwellian. The challenge arises from the nonlocal quadratic interaction term defined in the unbounded domain of velocity space. Analyzing this term on an unbounded domain requires the inclusion of a truncation function, which demands delicate analysis techniques. As a generalization of this analysis, we also provide proof of the asymptotic preserving property when using micro-macro decomposition-based neural networks.

\end{abstract}
\maketitle

\section{Introduction}
	
	Classical numerical methods (e.g., finite difference/volume/element methods) succeed in solving low-dimensional PDEs but encounter challenges in both theoretical analysis and numerical implementation when applied to high-dimensional counterparts. The utilization of machine learning techniques, particularly deep neural networks (DNNs), for solving high-dimensional PDEs, has rapidly gained attention and shown promising results in various applications (e.g., \cite{E2018, raissi2019physics, lu2021learning, li2020fourier}). In these approaches, DNNs are employed to minimize loss functions, posing a high-dimensional non-convex optimization challenge. The selection of appropriate loss functions is crucial in the context of solving PDEs (e.g., \cite{E2018, raissi2019physics, li2020fourier, beck2020overview, liao2019deep, deepGalerkin2018, zang2020weak, cai2021least, lyu2020mim}). For comprehensive reviews and references, we refer to, e.g., \cite{HJJL, CLM, LY, DHY, lou2021physics}.
	
	Kinetic equations serve as bridges between continuum and atomistic models \cite{Cerci}. The major challenge lies in the curse of dimensionality, as kinetic equations describe the evolution of probability density functions of a large number of particles and are defined in phase space, typically a six-dimensional problem plus the time dimension. When uncertainties are considered, the dimensionality can be even higher (e.g., \cite{Jin-Pareschi-Book, HuJin-UQ, poette2019gpc, poette2022numerical}). The DNN approach offers several advantages. Firstly, it can effectively handle high-dimensional PDEs due to its robustness and expressiveness. Secondly, it operates as a mesh-free method, enabling seamless navigation of complex domains and geometries. Thirdly, it offers user-friendly implementation by treating the residual error of PDEs as the loss function, eliminating the need to construct numerical schemes for approximating derivatives via automatic differentiation techniques. Therefore, DNN-based methods are ideal choices to tackle these challenges.
	
	Recently, \cite{jin2023asymptotic} introduced a neural network-based method for computing time-dependent linear transport equations with diffusive scaling and uncertainties. The goal of the network was to address the computational challenges associated with the curse of dimensionality and multiple scales of the problem. Despite the success of the numerical methods, rigorous analysis is still lacking, which motivates our studies in this paper. Since the analysis tools that shall be developed will generally work for Boltzmann-type equations, we do not limit our studies only to linear transport equations. Instead, we consider the classical Boltzmann equation and provide a rigorous error analysis for physics-informed neural networks (PINNs) in approximating it near a global Maxwellian. The challenge arises from the nonlocal quadratic interaction term defined in the unbounded domain of velocity space.  Additionally, we will provide proof of the asymptotic preserving property when using micro-macro decomposition-based neural networks.
	
	\subsection{Preliminaries} The classical Boltzmann equation takes the form:
	%\be
	\begin{equation}\la{eq:boltzmann}
	\partial_t f + \xi\cdot\nabla_x f=Q(f,f), \quad\text{with } (t,x,\xi)\in\R_+\times\R^n\times\R^n.
	\end{equation}
	%\ee
	Here, $f=f(t,x,\xi)$ is the mass density function of gas particles at time $t$ with position $x=(x_1,\cdots,x_n)$ and velocity $\xi=(\xi_1,\cdots,\xi_n)$. The {collision operator} $Q$ is defined by
	%\be
	\begin{equation}\la{eq:collision}
	Q(f,g)=\frac12 \int_{\R^n\times S^{n-1}} q(v,\beta) (  f'g'_* + f'_* g' - fg_* - f_* g ) \,d\xi_* d\omega,
	\end{equation}
	%\ee
	where 
	\begin{equation}
	f=f(t,x,\xi),\quad f'=f(t,x,\xi'), \quad f'_*=f(t,x,\xi'_*), \quad f_*=f(t,x,\xi_*),
	\end{equation}
	and similar for $g$, with the after collision velocities 
	\begin{equation}
	\xi'=\xi - \left(v\cos\beta\right)\omega, \quad\text{and}\quad 
	\xi'_*=\xi_* + \left(v\cos\beta\right)\omega,
	%\xi'=\xi - \left((\xi-\xi_*)\cdot\omega\right)\omega, \quad\text{and }
	%\xi'_*=\xi_* + \left((\xi-\xi_*)\cdot\omega\right)\omega,
	\end{equation}
	where $\omega\in S^{n-1}$, 
    $v=|\xi-\xi_*|$, and  $\cos\beta=\frac1v(\xi-\xi_*)\cdot\omega$.
	% \begin{equation}
	% v=|\xi-\xi_*|, \quad\text{and } \cos\theta=\frac1v(\xi-\xi_*)\cdot\omega.
	% \end{equation}
	Notice that the momentum and energy are conserved through the collision, i.e.,
	\begin{equation}
	\xi+\xi_*=\xi'+\xi'_* ,
	\quad\text{and }
	\abs{\xi}^2+\abs{\xi_*}^2=\abs{\xi'}^2+\abs{\xi'_*}^2 .
	\end{equation}
	
	%%% the properties of the collision operator
	It is well known that the collision operator satisfies the following properties:
	\begin{prop} 
		The collision operator satisfies the following properties:
		\begin{enumerate}
			\item It is held that 
			% \begin{equation}
			% \varphi_0(\xi)=1,\quad \varphi_i(\xi)=\xi_i\;(i=1,2,\cdots,n),\quad \varphi_{n+1}=\frac12\abs{\xi}^2,
			% \end{equation}
			% one has 
			\begin{equation}
			\int_{\R^n}\left(\begin{array}{c}1\\\xi\\\frac12\abs{\xi}^2\end{array}\right)\,Q(f,f)\,d\xi\equiv0.
            %\quad\text{for all }j=0, 1, 2, \cdots, n, n+1.
			\end{equation}
			\item The kernel of $Q$ is spanned by the \emph{Maxwellian}, 
			\begin{equation}
			M(\xi):=\left(\frac{\rho}{2\pi T}\right)^{d/2}\exp\left(-\frac{\abs{\xi-u}^2}{2T}\right).
			\end{equation}
			\item Entropy dissipation:
			\begin{equation}
			\int_{\R^n}\log{f}\,Q(f,f)\,d\xi\leq0.
			\end{equation}
		\end{enumerate}
	\end{prop}

	%%% Linearized boltzmann equation
	If we are looking for the solution $f$ of the Boltzmann equation \eqref{eq:boltzmann} near a global Maxwellian $M={(2\pi)^{-n/2}}\exp(-\abs{\xi}^2/2)$, then the solution can be written as
	\begin{equation}
	f = M+M^{1/2}u,
	\end{equation}
	from which one can obtain the following Boltzmann equation for $u$
	\be\la{eq:linear_boltzmann} 
	\pa_t u + \xi \cdot \na_x u + \nu(\xi) u = K(u) + \Gamma (u,u),
    %\pa_t u + \xi \cdot \na_x u = L(u) + \Gamma (u,u),
	\ee
	where
	\begin{align}
%	L(u)&:=
 -\nu(\xi)u+Ku:=2M^{-1/2}Q(M,M^{1/2}u), \;\;\text{and}\;\;	\Gamma(u,u):=M^{-1/2}Q(M^{1/2}u,M^{1/2}u).
	\end{align}

	We consider the Boltzmann equations \eqref{eq:linear_boltzmann}  
	on $[0,\infty) \times \TT^3 \times \R^3$.
    Under the angular cut-off assumption in \cite{ukai2007}, the operator $K: L^2 \rightarrow L^2 \cap L^{\infty}$ is linear and bounded, and $\nu$ obeys 
	\be 
	\nu_0 (1+ |\xi|)^{\gamma} \le \nu(\xi) \le \nu_1 (1+|\xi|)^{\gamma},  
	\ee for some positive constants $\nu_0$ and $\nu_1$ and some $\gamma \in [0,1]$. For a fixed $R>0$, we denote by $Q_R$ the cube in $\R^3$ centered at 0 with diameter $R$. For a set $A \subset \R^3$, we denote by $\chi_{A}$ the characteristic function of $A$.

	\subsection{PINNs for Boltzmann equation}
	In this paper, we estimate the errors resulting from approximations of solutions to the Boltzmann equations \eqref{eq:linear_boltzmann} by physics-informed neural networks (PINNs) \cite{raissi2019physics, jin2023asymptotic}. Let $u_\theta$ be a neural network approximation to the solution of the Boltzmann equation \eqref{eq:linear_boltzmann}, with $\theta$ denoting the trainable parameters of the network.
	
	To this end, we fix a sufficiently large diameter $R> 0$ and consider the following different types of residuals: the pointwise \emph{PDE residual}  
	\begin{equation*}
	\mathcal{R}_{i}[R; \theta] (t,x, \xi) = \pa_t u_{\theta} 
	+ \xi \cdot \na_{x} u_{\theta}
	+ \nu(\xi) u_{\theta}
	- K(u_{\theta} \chi_{Q_R})
	- \Gamma (u_{\theta} \chi_{Q_R}, u_{\theta} \chi_{Q_R}),
	\end{equation*}
	the \emph{initial residual} 
	\begin{equation*}
	\mathcal{R}_{t}[\theta] (x, \xi) = u_{\theta}(0, x, \xi) - u (0, x, \xi),
	\end{equation*} and the \emph{boundary residual} 
	\begin{align*}
	\mathcal{R}_{b}[\theta](t,x, \xi) 
	= \mathcal{R}_{b}[\theta](t,x_1,x_2, x_3, \xi) 
	&= (u_{\theta} (t, x_1, x_2, \pi, \xi) - u_{\theta} (t, x_1, x_2, - \pi, \xi ))^2 
	\\&\quad+ (u_{\theta} (t, x_1, \pi, x_3, \xi) - u_{\theta} (t, x_1, -\pi, x_3, \xi))^2
	\\&\quad\quad+ (u_{\theta} (t, \pi, x_2, x_3, \xi) - u_{\theta} (t, -\pi, x_2, x_3, \xi ))^2, 
	\end{align*}
	where $x_1$, $x_2$ and $x_3$ are the three components of $x \in \TT^3$.
	
	We define two types of errors:  the \emph{generalization error} and the  \emph{total error}. For a fixed diameter $R > 0$, the generalization error is defined as
	\be \label{generaler}
	\mathcal{E}_{G}[R;\theta] = \left(\mathcal{E}_{G}^i [R;\theta]^2 + \mathcal{E}_{G}^t [R;\theta]^2 + \mathcal{E}_{G}^b [R;\theta]^2 +  \lambda_R \mathcal{E}_{G}^p [R;\theta]^2 \right)^{\fr{1}{2}},
	\ee 
	where 
	\begin{align*}
	& \mathcal{E}_{G}^i [R;\theta]^2 
	= \int_{0}^{T} \int_{\TT^3} \int_{Q_R} \mathcal{R}_{i}[R;\theta]^2 d\xi dx dt,
	\\
	&\mathcal{E}_{G}^t [R;\theta]^2
	= \int_{\TT^3} \int_{Q_R} \mathcal{R}_{t}[\theta]^2 d\xi dx,
	\\
	&\mathcal{E}_{G}^b [R;\theta]^2
	=  \int_0^T \int_{\pa \TT^3} \int_{Q_R} \mathcal{R}_b[\theta]^2 d\xi d\sigma(x) dt,
	\\
	&\mathcal{E}_{G}^p [R;\theta]^2
	= (1+R)^{\gamma} \int_{0}^{T} \|u_{\theta}(t, \cdot, \xi)\|_{H^{2} (\TT^3 \times Q_R)}^2  dt.
	\end{align*}
	Here, $\lambda_R$ is a constant to be determined later. 
	
	For a fixed $R > 0 $, the total error is defined by 
	\be 
	\mathcal{E}[R; \theta]^2 = \int_{0}^{T} \int_{\TT^3} \int_{Q_R} |u - u_{\theta}|^2 d\xi dx dt. \nonumber
	\ee 
	
	\subsection{Main results} Let $T>0$ be an arbitrary positive time. Let $u$ be a solution to the Boltzmann equation \eqref{eq:linear_boltzmann} obeying $ \|(1+|\xi|)^{\frac{\gamma}{2}}u \|_{L^{\infty} (0,T; L^{\infty}(\TT^3; L^2(\RR^3)))}^2 < \infty$, and $\widehat{u}$ be \texttt{tanh} neural networks (see, e.g., \cite{de2021approximation}) approximating $u$. Then one can have the estimate $\mathcal{E}[R; \theta]\leq C \mathcal{E}_G[R; \theta]$ with the bound $C$ specified later in Theorem~\ref{thm:total}. This implies that as the loss function converges to zero, the neural network approximation converges to the true solution of the Boltzmann equation \eqref{eq:linear_boltzmann} near a global Maxwellian. Moreover, we also show that (in Theorem~\ref{generalmain}), for any $\epsilon>0$, under proper conditions, there exists a \texttt{tanh} neural network $\widehat{u}$ such that the generalization error $\mathcal{E}_G[R;\theta] \le \epsilon.$  Additionally, as a generalization of this analysis, we consider the rescaled Boltzmann equation and prove the asymptotic preserving property, whose numerical study is presented in the APNNs \cite{jin2023asymptotic}. 
	
\subsection{Related literature} The studies of using NNs to accelerate the computation of the Boltzmann equations include, e.g., \cite{han_uniformly_2019} for the moments closure method, \cite{xiao_using_2021} for the data-driven method, and \cite{lou_physics-informed_2021,li_physics-informed_2022,li_physics-informed_2023,jin2023asymptotic} for PINNs.
Recent research has focused not only on evaluating the computational efficiency and accuracy of PINNs in numerically solving PDEs but also on comprehensively studying their theoretical guarantees. A series of recent studies \cite{de2021approximation, de2022generic, de2022error, de2024error, mishra2022estimates, mishra2023estimates, jiao2021rate, jiao2023rate, jiao2023improved, guo2022mc, pang2020npinn, biswas2022error, hu2023higher, abdo2024accuracy} have thoroughly examined the error analysis of PINNs, investigating their capability to approximate various types of PDEs. These works shed light on our current paper by providing valuable insights into the accuracy and performance of PINNs across different problem domains. In this paper, we focus on the analysis of PINNs for the Boltzmann equation on an unbounded domain, which has not been investigated in the literature.

\subsection{Organization of the paper} In Section~\ref{sec:nonlocal}, as a preparation, we provide the estimates for the nonlocal operators $K$ and $\Gamma$; we estimate the generalization error of PINNs in Section~\ref{sec:gen_err} and total error in Section~\ref{sec:totalerr}. We consider the multiscale Boltzmann-type equation and prove the asymptotic preserving property in Section~\ref{sec:AP}. In the end, we make conclusive remarks in Section~\ref{sec:conclusion}.

\section{Local Estimates for the Nonlocal Operators \texorpdfstring{$K$}{} and \texorpdfstring{$\Gamma$}{}}\label{sec:nonlocal}

The incorporation of a characteristic function into the PDE residual results in some error that needs to be quantified. Due to their nonlocality, the operators $K$ and $\Gamma$ require some technical investigation. In this section, we establish new quantitative estimates for these operators that are crucial to studying the generalization and total errors in the upcoming sections.

\beg{prop}[Local Estimates for $K$]
Let $R > 0$. Let $v(x,\xi) \in L^2(\TT^3 \times \R^3)$. Then it holds that 
\be \la{linear1}
\|Kv - K(v \chi_{Q_R})\|_{L^2(\TT^3 \times Q_R)} 
\le C\|v\|_{L^2(\TT^3 \times \R^3 \setminus Q_R)}.
\ee  Moreover, if $\tilde{v} \in L^2(\TT^d \times Q_R)$, then it holds that 
\be \la{linear2}
\|K(v \chi_{Q_R}) - K(\tilde{v}\chi_{Q_R})\|_{L^2(\TT^3\times Q_R)} \le C \|v - \tilde{v}\|_{L^2(\TT^3 \times Q_R)}.
\ee 
\end{prop} 

\beg{proof}. In view of the boundedness of the linear operator $K(x, \cdot)$ on $L^2(\RR^3)$, we have
\bes
\beg{aligned}
\|Kv - K(v \chi_{Q_R})\|_{L^2(\TT^3 \times Q_R)}^2 
&= \int_{\TT^3} \int_{Q_R} K(v - v\chi_{Q_R})^2 d\xi dx
\le  \int_{\TT^3} \int_{\R^3} K(v - v\chi_{Q_R})^2 d\xi dx
\\&\le C  \int_{\TT^3} \int_{\R^3} (v - v\chi_{Q_R})^2 d\xi dx
= C\int_{\TT^3} \int_{\RR^3 \setminus Q_R} v^2 d\xi dx,
\end{aligned}
\ees which gives \eqref{linear1}. Similarly, we estimate 
\bes
\beg{aligned}
\|K(v \chi_{Q_R}) - K(\tilde{v}\chi_{Q_R})\|_{L^2(\TT^3\times Q_R)}^2
&= \int_{\TT^3} \int_{Q_R} K((v-\tilde{v})\chi_{Q_R})^2 d\xi dx
\le \int_{\TT^3} \int_{\R^3} K((v-\tilde{v})\chi_{Q_R})^2 d\xi dx
\\&\le C\int_{\TT^3} \int_{\R^3} (v-\tilde{v})^2\chi_{Q_R}^2 d\xi dx
= C\int_{\TT^3} \int_{Q_R} (v- \tilde{v})^2 d\xi dx,
\end{aligned}
\ees and we deduce \eqref{linear2}.
\end{proof}

% p=\infty, \alpha=1 case, Grad 1965 asymptotic equivalence of the NS and nonlinear Boltzmann
% p=2, \alpha=1/2, Golse Perthame Sulem 1987
\beg{lem}[{\cite[Theorem 1.2.3]{ukai2007}}] \label{wholespace} For $p \in [1,\infty], \alpha \in [0,1]$, there is a constant $C_0>0$ such that for any $v, \tilde{v} \in L^p_{\xi}(\R^3)$, it holds that 
\be 
\|\nu^{-\alpha} \Gamma(v, \tilde{v})\|_{L^p_{\xi}(\R^3)} 
\le C_0 \left(\|\nu^{1-\alpha} v\|_{L_{\xi}^p(\R^3)} \|\tilde{v}\|_{L_{\xi}^p(\R^3)} + \|\nu^{1-\alpha} \tilde{v}\|_{L_{\xi}^p(\R^3)} \|v\|_{L_{\xi}^p(\R^3)}\right).
\ee 
\end{lem}

\beg{prop}[Local Estimates for $\Gamma$] 
Let $R>0$. Let $\alpha \in [0,1]$.
\begin{enumerate}
\item[(i)] Suppose $v(x, \xi)$ obeys \be
(1+|\xi|)^{(1-\alpha)\gamma} v(x, \xi) \in L^{\infty}(\TT^3 ; L^2(\R^3)),
\ee and 
\be 
\int_{\TT^3} \int_{\R^3 \setminus Q_R} (1+ |\xi|)^{2\gamma(1-\alpha)} |v|^2 d\xi dx < \infty.
\ee Then it holds that 
\be \la{bilinear1}
\beg{aligned}
&\|\nu^{-\alpha} \Gamma (v - v\chi_{Q_R}, v)\|_{L^2(\TT^3 \times Q_R)}
\\&\quad\le C \|(1+|\xi|)^{(1-\alpha)\gamma} v\|_{L^{\infty}(\TT^3; L^2(\R^3))} \|(1+|\xi|)^{(1-\alpha)\gamma} v\|_{L^2(\TT^3 \times \R^3 \setminus Q_R)},
\end{aligned}
\ee and 
\be \la{bilinear2}
\beg{aligned}
&\|\nu^{-\alpha} \Gamma (v\chi_{Q_R}, v - v\chi_{Q_R})\|_{L^2(\TT^3 \times Q_R)}
\\&\quad\le C \|(1+|\xi|)^{(1-\alpha)\gamma} v\|_{L^{\infty}(\TT^3; L^2(\R^3))} \|(1+|\xi|)^{(1-\alpha)\gamma}v\|_{L^2(\TT^3 \times \R^3 \setminus Q_R)}.
\end{aligned}
\ee
\item[(ii)] Suppose $v(x,\xi) \in L^2(\TT^3 \times Q_R)$  such that 
\be 
(1+|\xi|)^{(1-\alpha)\gamma} v(x,\xi) \in L^{\infty}(\TT^3 ; L^2(\RR^3)),
\ee and $\tilde{v}(x,\xi) \in L^2(\TT^3 \times Q_R) \cap L^{\infty}(\TT^3;L^2(Q_R))$. Then it holds that 
\begin{multline}\la{bilinear3}
    \|\nu^{-\alpha} \Gamma ((v-\tilde{v})\chi_{Q_R}, v \chi_{Q_R}) \|_{L^2(\TT^3 \times Q_R)}^2
\\\le C \left(\|(1+|\xi|)^{(1-\alpha)\gamma}v\|_{L^{\infty}(\TT^3; L^2(Q_R))}^2 + (1+R)^{2\gamma(1-\alpha)} \|\tilde{v}\|_{L^{\infty}(\TT^3; L^2(Q_R))}^2 \right)\\{}\times\|v- \tilde{v}\|_{L^2(\TT^3 \times Q_R)}^2,
\end{multline}
and 
\begin{multline} \la{bilinear4}
\|\nu^{-\alpha} \Gamma (\tilde{v}\chi_{Q_R},(v-\tilde{v})\chi_{Q_R}) \|_{L^2(\TT^3 \times Q_R)}^2
\\\le C \left(\|(1+|\xi|)^{(1-\alpha)\gamma}v\|_{L^{\infty}(\TT^3;L^2(Q_R))}^2 + (1+R)^{2\gamma(1-\alpha)} \|\tilde{v}\|_{L^{\infty}(\TT^3; L^2(Q_R))}^2 \right) \\{}\times\|v- \tilde{v}\|_{L^2(\TT^3 \times Q_R)}^2.
\end{multline}
\end{enumerate}
\end{prop}

\beg{proof}. In view of Lemma \ref{wholespace} and the estimate $\nu(\xi) \le (1+|\xi|)^{\gamma}$ that holds for all $\xi \in \R^3$, we have 
\bes 
\beg{aligned}
&\|\nu^{-\alpha} \Gamma (v - v\chi_{Q_R}, v)\|_{L^2(\TT^3 \times Q_R)}^2
= \int_{\TT^3} \int_{Q_R} |\nu^{-\alpha} \Gamma (v - v\chi_{Q_R}, v) |^2 d\xi dx
\\&\quad\le \int_{\TT^3} \int_{\R^3} |\nu^{-\alpha} \Gamma (v - v\chi_{Q_R}, v) |^2 d\xi dx
\\&\quad\le C\int_{\TT^3} \left(\|\nu^{1-\alpha} v (1-\chi_{Q_R}) \|_{L^2(\R^3)}^2 \|v\|_{L^2(\R^3)}^2 + \| v (1-\chi_{Q_R}) \|_{L^2(\R^3)}^2 \|\nu^{1-\alpha} v\|_{L^2(\R^3)}^2\right)dx
\\&\quad= C\int_{\TT^3} \left(\|\nu^{1-\alpha} v  \|_{L^2(\R^3 \setminus Q_R)}^2 \|v\|_{L^2(\R^3)}^2 + \| v  \|_{L^2(\R^3 \setminus Q_R)}^2 \|\nu^{1-\alpha} v\|_{L^2(\R^3)}^2\right)dx 
\\&\quad\le C\int_{\TT^3} \|(1+|\xi|)^{(1-\alpha)\gamma}v\|_{L^2(\R^3 \setminus Q_R)}^2 \|(1+|\xi|)^{(1-\alpha)\gamma} v \|_{L^2(\R^3)}^2
\\&\quad\le C\|(1+|\xi|)^{(1-\alpha)\gamma} v  \|_{L^{\infty}(\TT^3 ; L^2(\RR^3))}^2 \|(1+|\xi|)^{(1-\alpha)\gamma} v \|_{L^2(\TT^3 \times \R^3 \setminus Q_R)}^2,
\end{aligned}
\ees yielding \eqref{bilinear1}. The proof of \eqref{bilinear2} is similar and will be omitted. Now we prove the bilinear estimate \eqref{bilinear3}. Indeed, we have
\be 
\beg{aligned}
&\|\nu^{-\alpha} \Gamma ((v-\tilde{v})\chi_{Q_R}, v \chi_{Q_R}) \|_{L^2(\TT^3 \times Q_R)}^2
= \int_{\TT^3} \int_{Q_R} |\nu^{-\alpha} \Gamma((v-\tilde{v})\chi_{Q_R}, v\chi_{Q_R}) |^2 d\xi dx
\\&\quad\le \int_{\TT^3} \int_{\R^3} |\nu^{-\alpha} \Gamma((v-\tilde{v})\chi_{Q_R}, v\chi_{Q_R}) |^2 d\xi dx
\\&\quad\le C\int_{\TT^3} \left(\|\nu^{1-\alpha} (v -\tilde{v})\chi_{Q_R} \|_{L^2(\R^3)}^2 \|v \chi_{Q_R} \|_{L^2(\R^3)}^2 +  \|\nu^{1-\alpha} v \chi_{Q_R} \|_{L^2(\R^3)}^2 \|(v-\tilde{v})\chi_{Q_R} \|_{L^2(\R^3)}^2 \right) dx
\\&\quad= C\int_{\TT^3} \left(\|\nu^{1-\alpha} (v -\tilde{v})\|_{L^2(Q_R)}^2 \|v  \|_{L^2(Q_R)}^2 +  \|\nu^{1-\alpha} v  \|_{L^2(Q_R)}^2 \|v-\tilde{v} \|_{L^2(Q_R)}^2 \right) dx
\\&\quad\le  C\int_{\TT^3} \|\nu^{1-\alpha} (v -\tilde{v})\|_{L^2(Q_R)}^2 \|v  \|_{L^2(Q_R)}^2 dx + C\|(1+|\xi|)^{(1-\alpha)\gamma} v\|_{L^{\infty}(\TT^3; L^2(Q_R))}^2 \|v - \tilde{v}\|_{L^2(\TT^3 \times Q_R)}^2.
\end{aligned}
\ee Several applications of the triangle inequality give rise to the following estimate,
\bes 
\beg{aligned}
&\int_{\TT^3} \|\nu^{1-\alpha} (v -\tilde{v})\|_{L^2(Q_R)}^2 \|v  \|_{L^2(Q_R)}^2 dx 
\\&\quad\le C\int_{\TT^3} \|\nu^{1-\alpha} (v -\tilde{v})\|_{L^2(Q_R)}^2 \|v -\tilde{v} \|_{L^2(Q_R)}^2 dx 
+ C\int_{\TT^3} \|\nu^{1-\alpha} (v -\tilde{v})\|_{L^2(Q_R)}^2 \|\tilde{v}  \|_{L^2(Q_R)}^2 dx 
\\&\quad\le C\int_{\TT^3} \left(\|\nu^{1-\alpha} v \|_{L^2(Q_R)}^2 + \|\nu^{1-\alpha} \tilde{v}\|_{L^2(Q_R)}^2\right)\|v -\tilde{v} \|_{L^2(Q_R)}^2 dx 
\\&\quad\quad\quad\quad+ C(1+R)^{2\gamma(1-\alpha)} \int_{\TT^3} \|v - \tilde{v}\|_{L^2(Q_R)}^2 \|\tilde{v}\|_{L^2(Q_R)}^2 dx
\\&\quad\le C \|(1+|\xi|)^{(1-\alpha)\gamma} v\|_{L^{\infty}(\TT^3;L^2(Q_R))}^2 \|v - \tilde{v}\|_{L^2(\TT^3 \times Q_R)}^2 
\\&\quad\quad\quad\quad+ (1+R)^{2\gamma(1-\alpha)} \|\tilde{v}\|_{L^{\infty}(\TT^3; L^2(Q_R))}^2 \|v- \tilde{v}\|_{L^2(\TT^3 \times Q_R)}^2.
\end{aligned}
\ees As a consequence, we infer that \eqref{bilinear3} holds. The proof of \eqref{bilinear4} is similar and will be omitted. 
\end{proof}

\section{Generalization Error Estimates}\label{sec:gen_err}

\beg{lem}[{\cite[Lemma A.1]{hu2023higher}}] Let $\Omega = \prod\limits_{i=1}^{d} [a_i, b_i]$. Suppose $f \in H^m(\Omega).$ Let $N>5$ be an integer. Then there exists a \texttt{tanh} neural network $\widehat{f}^N$, such that  for any $k \in \left\{0, 1, \dots, m-1 \right\}$,
$$\|f - \widehat{f}^N\|_{H^k(\Omega)} \le C_{k,d,f,\Omega} (1 + \ln^k N)N^{-m+k}.$$
\end{lem}

\beg{Thm} \label{generalmain} Suppose the solution $u$ to the Boltzmann equation obeys $\|(1+|\xi|)^{\gamma} u \|_{H^4([0,T] \times \TT^3 \times \RR^3)} <\infty$. Let $\epsilon >0$. There exists a real number $\lambda_R >0$, an integer $N>5$, and a \texttt{tanh} neural network $\widehat{u}= u_{\theta}^N$ 
such that $\mathcal{E}_G[R;\theta] \le \epsilon.$ 
\end{Thm}

\begin{proof}. As the solution $u$ to the Boltzmann equation \eqref{eq:linear_boltzmann} obeys 
\be
\int_{0}^{T}\int_{\TT^3} \int_{\RR^3} (1+|\xi|)^{2\gamma} |u|^2 d\xi dx dt < \infty, 
\ee it follows from the Lebesgue Dominated Convergence Theorem that there exists $R_0 > 0$ such that 
\be \label{longr}
\int_{0}^{T} \int_{\TT^3} \int_{\RR^3 \setminus Q_R} (1+|\xi|)^{2\gamma} |u|^2 d\xi dx dt < \delta(\epsilon)
\ee for any $R \ge R_0$, where $\delta(\epsilon)>0$ is a positive constant depending only on $\epsilon$ that will be chosen later.

Since $u \in H^4([0,T] \times \TT^3 \times Q_R)$, there exists a \texttt{tanh} neural network  $\widehat{u}$ such that 
\be \label{truesolutionapprox}
\|u - \widehat{u}\|_{H^3([0,T] \times \TT^3 \times Q_R)} \le C_{T, R, u} (1+\ln^3 N) N^{-1}.
\ee  
Fix $R \ge R_0$. Let $\Omega = [0,T] \times \TT^3 \times Q_R$. 

\smallskip
\noindent{\bf{Step 1. Estimation of $\mathcal{E}_{G}^i [R;\theta]$}.} 
As $u$ solves \eqref{eq:linear_boltzmann}, we can rewrite the PDE residual as 
\be 
\beg{aligned}
\mathcal{R}_i[R;\theta]
&= (\pa_t \widehat{u} - \pa_t u)
+ \xi \cdot \na_x (\widehat{u} - u)
+ \nu(\xi)(\widehat{u} - u)
+ K (u - u \chi_{Q_R})
+ K((u - \widehat{u})\chi_{Q_R})
\\&\quad\quad+ \Gamma(u, u) - \Gamma (u \chi_{Q_R}, u \chi_{Q_R})
+ \Gamma (u \chi_{Q_R}, u \chi_{Q_R}) - \Gamma (\widehat{u} \chi_{Q_R}, \widehat{u}\chi_{Q_R}).
\end{aligned}
\ee We have 
\be 
\beg{aligned}
&\|\pa_t (\widehat{u} - u)\|_{L^2(\Omega)} + \|\xi \cdot \na_x (\widehat{u} - u)\|_{L^2(\Omega)} 
+ \|\nu (\xi) (\widehat{u} - u)\|_{L^2(\Omega)}
\\&\le \|\widehat{u} - u\|_{H^1(\Omega)} + R\|\widehat{u} - u\|_{H^1(\Omega)} + \nu_1 (1+R)^{\gamma} \|\widehat{u} - u\|_{L^2(\Omega)}
\le C(1+ R) (1+\ln^3 N) N^{-1}
\end{aligned}
\ee in view of the approximation \eqref{truesolutionapprox}. By making use of \eqref{linear1}, \eqref{linear2}, \eqref{longr}, and \eqref{truesolutionapprox}, we estimate the linear terms as follows,
\beg{align}
\|Ku - K(u \chi_{Q_R}) \|_{L^2(\Omega)}
&\le C\|u\|_{L^2([0,T] \times \TT^3 \times \R^3 \setminus Q_R)}
\le C\delta(\epsilon),
\\ \|K(u \chi_{Q_R}) - K(\widehat{u} \chi_{Q_R}) \|_{L^2(\Omega)} &\le C\|u - \widehat{u}\|_{L^2(\Omega)} 
\le C(1+\ln^3 N) N^{-1}.
\end{align}
Applying the local bilinear bounds \eqref{bilinear1} and \eqref{bilinear2} with $\alpha = 0$, we obtain  
\be 
\beg{aligned}
&\|\Gamma(u,u) - \Gamma (u \chi_{Q_R}, u \chi_{Q_R})\|_{L^2(\Omega)}
\\&\quad\quad\le \|\Gamma (u(1-\chi_{Q_R}), u)\|_{L^2(\Omega)} + \|\Gamma (u \chi_{Q_R}, u (1 - \chi_{Q_R}))\|_{L^2(\Omega)}
\\&\quad\quad\le C\|(1+|\xi|)^{\gamma} u\|_{L^{\infty}(0,T;L^{\infty}(\TT^3; L^2(\RR^3)))} \| (1+|\xi|)^{\gamma} u\|_{L^2([0,T] \times \TT^2 \times \RR^3 \setminus Q_R)}
\\&\quad\quad\le C\|(1+|\xi|)^{\gamma} u\|_{L^2(\RR^3; H^3([0,T] \times \TT^3)}
\delta(\epsilon)
\\&\quad\quad\le C\|(1+|\xi|)^{\gamma} u\|_{H^3([0,T] \times \TT^3 \times \RR^3)} \delta(\epsilon).
\end{aligned}
\ee Here we used the continuous embedding of $L^2(\RR^3; H^3([0,T] \times \TT^3)$ into $L^{\infty} (0,T; L^{\infty}(\TT^3; L^2(\RR^3)))$ to bound $\|(1+|\xi|)^{\gamma} u\|_{L^{\infty} (0,T; L^{\infty}(\TT^3; L^2(\RR^3)))}$ by a constant multiple of $\|(1+|\xi|)^{\gamma} u\|_{L^2(\RR^3; H^3([0,T] \times \TT^3)}$. Indeed, this latter fact follows from 
\be 
\beg{aligned}
\sup\limits_{[0,T] \times \TT^3} \int_{\RR^3} (1+|\xi|)^{2\gamma} |u|^2 d\xi 
&\le \int_{\RR^3} \sup\limits_{[0,T] \times \TT^3} (1+|\xi|)^{2\gamma} |u|^2  d\xi 
\\&\le C\int_{\RR^3} \|(1+|\xi|)^{\gamma}  u\|_{H^3([0,T] \times \TT^3)}^2 d\xi
\\&\le C\|(1+|\xi|)^{\gamma} 
 u\|_{L^2(\RR^3; H^3([0,T] \times \TT^3)}^2
 \end{aligned}
\ee that holds due to the four-dimensional continuous embedding of $H^3([0,T] \times\TT^3)$ into $L^{\infty}([0,T] \times \TT^3)$. Using \eqref{bilinear3} and \eqref{bilinear4}, we estimate 
\be 
\beg{aligned}
&\|\Gamma (u \chi_{Q_R}, u \chi_{Q_R}) - \Gamma (\widehat{u} \chi_{Q_R}, \widehat{u} \chi_{Q_R})\|_{L^2(\Omega)}
\\&\quad\quad\le \|\Gamma ((u - \widehat{u})\chi_{Q_R}, u \chi_{Q_R})\|_{L^2(\Omega)}
+ \|\Gamma (\widehat{u} \chi_{Q_R}, (u - \widehat{u})\chi_{Q_R})\|_{L^2(\Omega)}
\\&\quad\quad\le C\left(\| (1+|\xi|)^{\gamma} u\|_{L^{\infty}(0,T; L^{\infty}(\TT^3; L^2(Q_R)))} + (1+R)^{\gamma} \|\widehat{u} \|_{L^{\infty}(0,T; L^{\infty}(\TT^3; L^2(Q_R)))} \right) \|u - \widehat{u}\|_{L^2(\Omega)}
\\&\quad\quad\le C(1+R)^{\gamma} \left(\|u\|_{H^3(\Omega)} + \|\widehat{u}\|_{H^3(\Omega)} \right) \|u - \widehat{u}\|_{L^2(\Omega)}
\\&\quad\quad\le C(1+R)^{\gamma} \|u\|_{H^3(\Omega)}\|u - \widehat{u}\|_{L^2(\Omega)}
+ C(1+R)^{\gamma} \|\widehat{u} - u\|_{H^3(\Omega)}  \|u - \widehat{u}\|_{L^2(\Omega)}
\\&\quad\quad\le C(1+R)^{\gamma} (\|u\|_{H^3([0,T]\times \TT^3 \times \RR^3)} +1) \left(1+ \ln^3 N \right)^2 N^{-1}.
\end{aligned}
\ee Therefore, we infer that 
\be 
\mathcal{E}_G^i [R; \theta] \le (1+R) (\|u\|_{H^3([0,T]\times \TT^3 \times \RR^3)} +1) (1 + \ln^6 N) N^{-1}  + C (\|(1+|\xi|)^{\gamma} u\|_{H^3([0,T]\times \TT^3 \times \RR^3)} +1)\delta(\epsilon).
\ee

\smallskip
\noindent{\bf{Step 2. Estimation of $\mathcal{E}_{G}^t[R;\theta]$}.} By the trace theorem, we have
\be \nonumber
\beg{aligned}
\mathcal{E}_{G}^t[\ell;\theta]
&= \|\widehat{u}(0, x, \xi) - u (0, x, \xi)\|_{L^{2}(\TT^3 \times Q_R)}
\le \|\widehat{u} - u\|_{L^2(\pa \Omega)} \le C\|\widehat{u} - u\|_{H^1(\Omega)} 
\le C(1+ \ln^3 N) N^{-1}.
\end{aligned}
\ee

\smallskip
\noindent{\bf{Step 3. Estimation of $\mathcal{E}_{G}^b [R;\theta]$}.} By making use of the spatial periodic boundary conditions obeyed by $u$ and applying the trace theorem, we deduce that 
\beg{align*}
&\mathcal{E}_{G}^b [R;\theta]
\le C\|\widehat{u} - u\|_{L^2(\pa \Omega)} \le C\|\widehat{u} - u\|_{H^1(\Omega)} \le C(1+\ln^3 N)N^{-1}.
\end{align*}

\smallskip
\noindent{\bf{Step 4. Estimation of $\mathcal{E}_{G}^p [R;\theta]$}.} By subtracting and adding $u$, we bound 
\be 
\beg{aligned}
\mathcal{E}_{G}^p [R;\theta] ^2
&\le (1+R)^{\gamma}\int_{0}^{T} \|\widehat{u} - u\|_{H^{2}(\TT^3 \times Q_R)}^2 dt
+ C(1+R)^{\gamma} \int_{0}^{T} \|u\|_{H^{2}(\TT^3 \times Q_R)}^2dt
\\&\le C(1+R)^{\gamma} \|\widehat{u} - u\|_{H^2(\Omega)}^2 + C(1+R)^{\gamma} \|u\|_{H^2(\Omega)}^2.
\end{aligned}
\ee Thus, we obtain 
\be \nonumber
\mathcal{E}_{G}^p [R; \theta]  
\le C(1+R)^{\frac{\gamma}{2}} (1+\ln^3N) N^{-1} + C(1+R)^{\frac{\gamma}{2}} \|u\|_{H^2(\Omega)}.
\ee 

\smallskip
\noindent{\bf{Step 5. Estimation of $\mathcal{E}_{G}[R;\theta]$}.} Combining the bounds derived in Steps 1--4, we infer that 
\be 
\mathcal{E}_G[R; \theta] \le C\lambda^{\fr{1}{2}}(1+R)^{\fr{\gamma}{2}}\|u\|_{H^2(\Omega)}  + C_{T, R, u} (1+ \ln^6 N)N^{-1}+ C (\|(1+|\xi|)^{\gamma} u\|_{H^3([0,T]\times \TT^3 \times \RR^3)} +1)\delta(\epsilon) .
\ee
We choose $\lambda$ and $\delta (\epsilon)$ sufficiently small and $N$ sufficiently large such that 
\beg{align}
C\lambda^{\fr{1}{2}}(1+R)^{\fr{\gamma}{2}}\|u\|_{H^2(\Omega)}  &= \fr{\epsilon}{3},
\\ C\left(\|(1+|\xi|)^{\gamma} u \|_{H^3([0,T] \times \TT^3 \times \RR^3)}+ 1 \right)\delta(\epsilon) &= \frac{\epsilon}{3},
\\C_{R,T,u}(1+\ln^6 N)N^{-1} &< \frac{\epsilon}{3}.
\end{align}
For these choices, the generalization error $\mathcal{E}_G[R;\theta] \le \epsilon.$
\end{proof}

Note that the assumption used in Theorem~\ref{generalmain} can be satisfied by, for example, the case studied in \cite{alexandre2011global}, as specified in the following remark.
\beg{rem} Suppose $q(v, \theta) = b(\cos \theta)$ where $b(\cos \theta)$ is of order $\theta^{-2-2s}$ as $\theta \rightarrow 0^+$ with $0 < s < \frac{1}{2}$. Suppose $u_0 (x, \xi)$ obeys $M + \sqrt{M} u_0 \ge 0$ and 
\be 
(1+ |\xi|)^{\ell} u_0 (x, \xi) \in H_x^k(\TT^3) \times H_{\xi}^k(\RR^3),
\ee for some $k \ge 3, \ell \ge 3$. There exists $\epsilon_0 >0$ such that if
\be 
\|(1+ |\xi|)^{{\ell}} u_0\|_{H^k_{x, \xi}(\TT^3 \times \R^3)} \le \epsilon_0,  
\ee then the Boltzmann equations \eqref{eq:linear_boltzmann} has a unique global solution $u$ obeying
\be 
\sup\limits_{t \ge 0}  \|(1+ |\xi|)^{\ell} u(t,x,\xi)\|_{H^k_{x, \xi}(\TT^3 \times \R^3)}  < \infty. 
\ee Moreover, it holds that $u \in H^k ([0,T] \times \TT^3 \times \RR^3)$. %{\color{red} The result is stated on $\RR^3 \times \RR^3$ but I think it works on $\TT^3 \times \RR^3$ as these two cases are usually similar}
\end{rem}

\beg{rem}
Let $\mathcal{S}_i, \mathcal{S}_t, \mathcal{S}_b$ be sets of points representing the midpoints of cubes partitioning $[0,T] \times \TT^3 \times Q_R$, $\TT^3 \times Q_R$ and $[0,T] \times \partial \TT^3 \times Q_R$.  
We define the training error $\mathcal{E}_T$ as follows:
\be 
\mathcal{E}_{T}[R; \theta; \mathcal{S}] = \left(\mathcal{E}_T^i[R;\theta;\mathcal{S}_i]^2 + \mathcal{E}_T^t[R;\theta;\mathcal{S}_t]^2
+ \mathcal{E}_T^b[R;\theta;\mathcal{S}_b]^2
+ \lambda_R \mathcal{E}_T^p[R;\theta;\mathcal{S}_i]^2
\right)^2 
\ee where 
\begin{align*}
\mathcal{E}_T^i[R;\theta;\mathcal{S}_i]^2 &= \sum\limits_{P_n \in \mathcal{S}_i} \frac{1}{|\mathcal{S}_i|} \mathcal{R}_i[R;\theta](P_n)^2,  \hspace{1cm} \mathcal{E}_T^t[R;\theta;\mathcal{S}_t]^2
= \sum\limits_{Q_n \in \mathcal{S}_t} \frac{1}{|\mathcal{S}_t|} \mathcal{R}_t[R;\theta](Q_n)^2,\\\\ \mathcal{E}_T^b[R;\theta;\mathcal{S}_b]^2
 &= \sum\limits_{R_n \in \mathcal{S}_b} \frac{1}{|\mathcal{S}_b|} \mathcal{R}_b[R;\theta](R_n)^2, \hspace{1cm} \mathcal{E}_T^p[R;\theta;\mathcal{S}_i]^2 =\sum\limits_{P_n \in \mathcal{S}_i} \sum\limits_{|\alpha| \le 2}\frac{1}{|\mathcal{S}_i|}D^{\alpha}u_{\theta}(P_n)^2.   
\end{align*} Using the midpoint rule, one can show that 
\be \label{gentot}
|\mathcal{E}_G[R;\theta]^2 - \mathcal{E}_T[R;\theta; \mathcal{S}]^2|
\le C_i|\mathcal{S}_i|^{-2/(d+1)}
+ C_t |\mathcal{S}_t|^{-2/d}
+C_b |\mathcal{S}_b|^{-2/d}
\ee where $d$ is the dimension of the cube $[0,T] \times \TT^3 \times Q_R$ (here $d=7$) and $C_i, C_t$ and $C_b$ are constants depending on the $C^2$ regularity of the corresponding residuals. Indeed, the $C^2(([0,T] \times \TT^3 \times Q_R)$ norm of $\mathcal{R}_i$ depends on the 
$C^3(([0,T] \times \TT^3 \times Q_R)$ norm of the neural network $u_{\theta}$ and the $C^2(([0,T] \times \TT^3 \times Q_R)$ norm of $\nu u_{\theta}, K(u_{\theta} \chi_{Q_R})$ and $\Gamma(u_{\theta} \chi_{Q_R}, u_{\theta} \chi_{Q_R})$. The estimate \eqref{gentot} shows the convergence of the training error to the generalization error as the cardinalities of the quadrature sets get very large. We note that the midpoint rule is computationally more efficient than the trapezoidal rule in higher dimensions and has been widely employed in the literature to analyze the training error arising from PINNs approximations (for instance, see \cite{hu2023higher} and references therein).
\end{rem}

\section{Total Error Estimates}\label{sec:totalerr}

\beg{Thm}\label{thm:total} Let $T>0$ be an arbitrary positive time. Let $u$ be a solution to the Boltzmann equations \eqref{eq:linear_boltzmann} obeying $ \|(1+|\xi|)^{\frac{\gamma}{2}}u \|_{L^{\infty} (0,T; L^{\infty}(\TT^3; L^2(\RR^3)))} < \infty$, and $\widehat{u}$ be \texttt{tanh} neural networks approximating $u$. Then there exists a positive constant $C$ depending only on $T$ and a diameter $R_0 > 0$ such that for any $R \ge R_0$, it holds that 
\be  \label{TEestimate}
\beg{aligned}
&\mathcal{E}[R; \theta]^2
\le  C\left(1 + \|(1+|\xi|)^{\frac{\gamma}{2}}u \|_{L^{\infty} (0,T; L^{\infty}(\TT^3; L^2(\RR^3)))}^2 \right) \mathcal{E}_G [R;\theta]^2 
\\&\quad\quad\times \exp \left\{\int_{0}^{T} C \left(1 + \|(1+|\xi|)^{\frac{\gamma}{2}}u \|_{L^{\infty}(\TT^3; L^2(\RR^3))}^2 + \frac{\mathcal{E}_G[R;\theta]^2}{\lambda_R} \right)dt \right\}.
\end{aligned}
\ee
\end{Thm}

\beg{proof}.
In view of the regularity criterion 
\be 
\int_{0}^{T} \int_{\TT^3} \int_{\RR^3} (1+|\xi|)^{\gamma} |u|^2 d\xi dx dt < \infty,
\ee obeyed by $u$ and the Lebesgue Dominated Convergence Theorem, there is a diameter $R_0 \ge 1$ such that for any $R \ge R_0$, it holds that 
\be 
\int_{0}^{T} \int_{\TT^3} \int_{\RR^3 \setminus Q_R} \left(1 + |\xi|\right)^{\gamma} |u|^2  d\xi dx dt \le \mathcal{E}_G^t [1;\theta]^2.
\ee Since $\mathcal{E}_G^t[R;\theta]$ is an increasing function in $R$, it follows that 
\be \label{decayest}
\int_{0}^{T} \int_{\TT^3} \int_{\RR^3 \setminus Q_R} \left(1 + |\xi|\right)^{\gamma} |u|^2  d\xi dx dt \le \mathcal{E}_G^t [R;\theta]^2 \le \mathcal{E}_G[R;\theta]^2,
\ee for all $R \ge R_0$. Now we fix an $R \ge R_0$. The difference $u - \widehat{u}$ evolves in time according to 
\be 
\beg{aligned}
&\pa_t (u - \widehat{u})
+ \xi \cdot \na_x (u - \widehat{u})
+ \nu(\xi) (u - \widehat{u})
\\&\quad\quad= K(u) - K(u \chi_{Q_R})
+ K(u \chi_{Q_R}) - K(\widehat{u} \chi_{Q_R})
\\&\quad\quad\quad\quad+ \Gamma (u,u) - \Gamma(u \chi_{Q_R}, u \chi_{Q_R}) + \Gamma (u \chi_{Q_R}, u \chi_{Q_R}) - \Gamma (\widehat{u} \chi_{Q_R}, \widehat{u} \chi_{Q_R}) + \mathcal{R}_{i}[R;\theta],
\end{aligned}
\ee which, after multiplying by $u - \widehat{u}$ and integrating over $\TT^3 \times Q_R$,  gives rise to the energy equation 
\be 
\beg{aligned}
\fr{1}{2} \frac{d}{dt} \|u - \widehat{u}\|_{L^2(\TT^3 \times Q_R)}^2 
+ \|\sqrt{\nu} (u - \widehat{u})\|_{L^2(\TT^3 \times Q_R)}^2
&= \mathcal{O}_1
+ \mathcal{O}_2
+ \mathcal{O}_3
+ \mathcal{O}_4
+ \mathcal{O}_5,
\end{aligned} 
\ee where 
\beg{align*}
\mathcal{O}_1 &= \int_{\TT^3} \int_{Q_R} \left(K(u) - K(u \chi_{Q_R})\right) (u - \widehat{u}) d\xi dx,
\\ \mathcal{O}_2 &= \int_{\TT^3} \int_{Q_R} (K(u \chi_{Q_R}) - K(\widehat{u} \chi_{Q_R}))(u - \widehat{u}) d\xi dx,
\\ \mathcal{O}_3 &= \int_{\TT^3} \int_{Q_R} (\Gamma (u,u) - \Gamma(u \chi_{Q_R}, u \chi_{Q_R}))(u - \widehat{u}) d\xi dx,
\\ \mathcal{O}_4 &= \int_{\TT^3} \int_{Q_R} (\Gamma (u \chi_{Q_R}, u \chi_{Q_R}) - \Gamma (\widehat{u} \chi_{Q_R}, \widehat{u} \chi_{Q_R}))(u - \widehat{u}) d\xi dx,
\\ \mathcal{O}_5&= \int_{\TT^3} \int_{Q_R} \mathcal{R}_{i}[R;\theta] (u - \widehat{u}) d\xi dx.
\end{align*}
By the Cauchy-Schwarz inequality and the local linear estimate \eqref{linear1},  we estimate $\mathcal{O}_1$ as follows,
\be 
|\mathcal{O}_1| 
\le \|Ku - K(u\chi_{Q_R})\|_{L^2(\TT^3 \times Q_R)} \|u - \widehat{u}\|_{L^2(\TT^3 \times Q_R)} 
\le C\|u\|_{L^2(\TT^3 \times \R^3 \setminus Q_R)}  \|u - \widehat{u}\|_{L^2(\TT^3 \times Q_R)}.
\ee Applying Young's inequality for products, we obtain the bound
\be 
|\mathcal{O}_1| \le \fr{1}{2}\|u - \widehat{u}\|_{L^2(\TT^3 \times Q_R)}^2 + C\|u\|_{L^2(\TT^3 \times \R^3 \setminus Q_R)}^2.
\ee By making use of  \eqref{linear2}, we estimate 
\be 
|\mathcal{O}_2| 
\le \|K(u \chi_{Q_R}) - K(\widehat{u}\chi_{Q_R})\|_{L^2(\TT^3 \times Q_R)} \|u - \widehat{u}\|_{L^2(\TT^3 \times Q_R)}
\le C\|u - \widehat{u}\|_{L^2(\TT^3 \times Q_R)}^2.
\ee In order to obtain good control of $\mathcal{O}_3$, we decompose it into the sum of two terms, $\mathcal{O}_{31}$ and $\mathcal{O}_{32}$, where 
\be 
\mathcal{O}_{31} = \int_{\TT^3} \int_{Q_R} \nu^{-\fr{1}{2}}\left(\Gamma(u,u) - \Gamma(u \chi_{Q_R}, u) \right) \nu^{\fr{1}{2}}(u - \widehat{u}) d\xi dx,
\ee and 
\be 
\mathcal{O}_{32} = \int_{\TT^3} \int_{Q_R} \nu^{-\fr{1}{2}}\left(\Gamma(u \chi_{Q_R}, u) - \Gamma(u \chi_{Q_R}, u \chi_{Q_R}) \right) \nu^{\fr{1}{2}}(u - \widehat{u}) d\xi dx.
\ee We have 
\be 
\beg{aligned}
|\mathcal{O}_{31}| 
&\le \|\nu^{-\fr{1}{2}} (\Gamma(u,u) - \Gamma (u \chi_{Q_R}, u)) \|_{L^2(\TT^3 \times Q_R)} \|\nu^{\fr{1}{2}} (u - \widehat{u}) \|_{L^2(\TT^3 \times Q_R)}
\\&\le C\|(1+|\xi|)^{\frac{\gamma}{2}}u \|_{L^{\infty}(\TT^3; L^2(\RR^3))
)} \|(1+|\xi|)^{\frac{\gamma}{2}}u\|_{L^2(\TT^3 \times \RR^3 \setminus Q_R)} \|\nu^{\fr{1}{2}} (u - \widehat{u}) \|_{L^2(\TT^3 \times Q_R)}
\\&\le \frac{1}{8} \|\sqrt{\nu} (u - \widehat{u})\|_{L^2(\TT^3 \times Q_R)}^2 + C\|(1+|\xi|)^{\frac{\gamma}{2}}u \|_{L^{\infty}(\TT^3; L^2(\RR^3))}^2 \|(1+|\xi|)^{\frac{\gamma}{2}}u\|_{L^2(\TT^3 \times \RR^3 \setminus Q_R)}^2,
\end{aligned} 
\ee and 
\be 
\beg{aligned}
|\mathcal{O}_{32}| 
&\le \|\nu^{-\fr{1}{2}} (\Gamma (u \chi_{Q_R}, u) - \Gamma(u\chi_{Q_R}, u \chi_{Q_R})) \|_{L^2(\TT^3 \times Q_R)} \|\nu^{\fr{1}{2}} (u - \widehat{u}) \|_{L^2(\TT^3 \times Q_R)}
\\&\le C\|(1+|\xi|)^{\frac{\gamma}{2}}u \|_{L^{\infty}(\TT^3; L^2(\RR^3))
)} \|(1+|\xi|)^{\frac{\gamma}{2}}u\|_{L^2(\TT^3 \times \RR^3 \setminus Q_R)} \|\nu^{\fr{1}{2}} (u - \widehat{u}) \|_{L^2(\TT^3 \times Q_R)}
\\&\le \frac{1}{8} \|\sqrt{\nu} (u - \widehat{u})\|_{L^2(\TT^3 \times Q_R)}^2 + C\|(1+|\xi|)^{\frac{\gamma}{2}}u \|_{L^{\infty}(\TT^3; L^2(\RR^3))}^2 \|(1+|\xi|)^{\frac{\gamma}{2}}u\|_{L^2(\TT^3 \times \RR^3 \setminus Q_R)}^2,
\end{aligned}
\ee as a consequence of the local bilinear estimates \eqref{bilinear1} and \eqref{bilinear2} and the bound \eqref{decayest}. As for the bilinear term $\mathcal{O}_4$ involving differences between the actual solution and neural networks, we also perform a decomposition $\mathcal{O}_4 = \mathcal{O}_{41} + \mathcal{O}_{42}$ where 
\be 
\mathcal{O}_{41} = \int_{\TT^3} \int_{Q_R} \nu^{-\fr{1}{2}}\left(\Gamma (u \chi_{Q_R}, u \chi_{Q_R}) - \Gamma (\widehat{u}\chi_{Q_R}, u \chi_{Q_R}) \right) \nu^{\fr{1}{2}}(u - \widehat{u}) d\xi dx,
\ee and 
\be 
\mathcal{O}_{42} = \int_{\TT^3} \int_{Q_R} \nu^{-\fr{1}{2}} (\Gamma(\widehat{u} \chi_{Q_R},u
\chi_{Q_R})- \Gamma(\widehat{u} 
\chi_{Q_R},\widehat{u} 
\chi_{Q_R})) \nu^{\fr{1}{2}} (u - \widehat{u}) d\xi dx.
\ee We bound 
\be 
\beg{aligned}
|\mathcal{O}_{41}| 
&\le \fr{1}{8} \|\sqrt{\nu} (u - \widehat{u})\|_{L^2(\TT^3 \times Q_R)}^2 + C\|\nu^{-\fr{1}{2}}\left(\Gamma (u \chi_{Q_R}, u \chi_{Q_R}) - \Gamma (\widehat{u}\chi_{Q_R}, u \chi_{Q_R}) \right)\|_{L^2(\TT^3 \times Q_R)}^2 
\\&\le \fr{1}{8} \|\sqrt{\nu} (u - \widehat{u})\|_{L^2(\TT^3 \times Q_R)}^2 + C \|(1+|\xi|)^{\frac{\gamma}{2}} u\|_{L^{\infty}(\TT^3; L^2(Q_R))}^2  \|u - \widehat{u}\|_{L^2(\TT^3 \times Q_R)}^2
\\&\quad\quad+ C (1+R)^{\gamma} \|\widehat{u}\|_{L^{\infty}(\TT^3; L^2(Q_R))}^2 \|u - \widehat{u}\|_{L^2(\TT^3 \times Q_R)}^2,
\end{aligned}
\ee and 
\be 
\beg{aligned}
|\mathcal{O}_{42}| 
&\le \fr{1}{8} \|\sqrt{\nu} (u - \widehat{u})\|_{L^2(\TT^3 \times Q_R)}^2 + \|\nu^{-\fr{1}{2}} (\Gamma(\widehat{u} \chi_{Q_R},u
\chi_{Q_R})- \Gamma(\widehat{u} 
\chi_{Q_R},\widehat{u} 
\chi_{Q_R}))\|_{L^2(\TT^3 \times Q_R)}^2 
\\&\le \fr{1}{8} \|\sqrt{\nu} (u - \widehat{u})\|_{L^2(\TT^3 \times Q_R)}^2 + C \|(1+|\xi|)^{\frac{\gamma}{2}} u\|_{L^{\infty}(\TT^3; L^2(Q_R))}^2  \|u - \widehat{u}\|_{L^2(\TT^3 \times Q_R)}^2
\\&\quad\quad+ C (1+R)^{\gamma} \|\widehat{u}\|_{L^{\infty}(\TT^3; L^2(Q_R))}^2 \|u - \widehat{u}\|_{L^2(\TT^3 \times Q_R)}^2,
\end{aligned}
\ee using the local bilinear estimates \eqref{bilinear3} and \eqref{bilinear4}. Finally, we have 
\be 
|\mathcal{O}_5| \le \fr{1}{2}\|u - \widehat{u}\|_{L^2(\TT^3 \times Q_R)}^2 + C\|\mathcal{R}_i[R;\theta] \|_{L^2(\TT^3 \times Q_R)}^2.
\ee 
Therefore, we obtain the differential inequality
\bes 
\beg{aligned}
&\fr{d}{dt} \|u - \widehat{u}\|_{L^2(\TT^3 \times Q_R)}^2 
+ \|\sqrt{\nu} (u - \widehat{u})\|_{L^2(\TT^3 \times Q_R)}^2
\\&\quad\quad\le C\left(1 + \|(1+|\xi|)^{\frac{\gamma}{2}}u \|_{L^{\infty}(\TT^3; L^2(\RR^3))}^2 \right)\|(1 + |\xi|)^{\fr{\gamma}{2}} u\|_{L^2(\TT^3 \times \RR^3 \setminus Q_R)}^2
\\&\quad\quad\quad+ C\left(1 + \|(1+|\xi|)^{\frac{\gamma}{2}}u \|_{L^{\infty}(\TT^3; L^2(\RR^3))}^2 + (1+R)^{\gamma} \|\widehat{u}\|_{L^{\infty}(\TT^3; L^2(Q_R))}^2  \right)\|u - \widehat{u}\|_{L^2(\TT^3 \times Q_R)}^2
\\&\quad\quad\quad\quad
+ C\|\mathcal{R}_i[R;\theta] \|_{L^2(\TT^3 \times Q_R)}^2.
\end{aligned}
\ees Gronwall's inequality yields the instantaneous bound
\bes 
\beg{aligned}
&\|(u - \widehat{u})(t)\|_{L^2(\TT^3 \times Q_R)}^2
\\&\quad\le \left[\mathcal{E}_G^t  [R;\theta]^2  + \mathcal{E}_G^i[R;\theta]^2 + \int_{0}^{T} \left(1 + \|(1+|\xi|)^{\frac{\gamma}{2}}u \|_{L^{\infty}(\TT^3; L^2(\RR^3))}^2 \right)\|(1 + |\xi|)^{\fr{\gamma}{2}} u\|_{L^2(\TT^3 \times \RR^3 \setminus Q_R)}^2 dt\right] 
\\&\quad\quad\quad\times \exp \left\{\int_{0}^{T} C \left(1 + \|(1+|\xi|)^{\frac{\gamma}{2}}u \|_{L^{\infty}(\TT^3; L^2(\RR^3))}^2 + (1+R)^{\gamma} \|\widehat{u}\|_{L^{\infty}(\TT^3; L^2(Q_R))}^2  \right)dt \right\}.
\end{aligned}
\ees 
Due to the three-dimensional continuous Sobolev embedding of $H^2(\TT^3)$ into $L^{\infty}(\TT^3)$, we have
\bes 
\sup\limits_{x \in \TT^3} \int_{Q_R} |\widehat{u}|^2 d\xi 
\le \int_{Q_R} \sup\limits_{x \in \TT^3} |\widehat{u}|^2 d\xi 
\le C\int_{Q_R} \|\widehat{u}\|_{H^2(\TT^3)}^2 d\xi 
\le C\|\widehat{u}\|_{H^2(\TT^3 \times Q_R)}^2,
\ees and consequently,
\bes 
(1+R)^{\gamma} \int_{0}^{T} \|\widehat{u}\|_{L^{\infty}(\TT^3; L^2(\Q_R))}^2 dt \le  \int_{0}^{T} C (1+R)^{\gamma}\|\widehat{u}\|_{H^2(\TT^3 \times Q_R)}^2 dt
= C\mathcal{E}_G^p [R;\theta]^2
\le \frac{C \mathcal{E}_{G}[R;\theta]^2}{\lambda_R}.
\ees Using \eqref{decayest}, we infer that 
\bes 
\beg{aligned}
&\|(u - \widehat{u})(t)\|_{L^2(\TT^3 \times Q_R)}^2
\le  C\left(1 + \|(1+|\xi|)^{\frac{\gamma}{2}}u \|_{L^{\infty} (0,T; L^{\infty}(\TT^3; L^2(\RR^3)))}^2 \right) \mathcal{E}_G [R;\theta]^2 
\\&\quad\quad\times \exp \left\{\int_{0}^{T} C \left(1 + \|(1+|\xi|)^{\frac{\gamma}{2}}u \|_{L^{\infty}(\TT^3; L^2(\RR^3))}^2 + \frac{\mathcal{E}_G[R;\theta]^2}{\lambda_R} \right)dt \right\}.
\end{aligned}
\ees Integrating in time from $0$ to $T$, we deduce that \eqref{TEestimate} holds.
\end{proof}

\section{Asymptotic preserving property}\label{sec:AP} 
This section addresses the asymptotic-preserving (AP) property of physics-informed neural networks (PINNs) in the multiscale regime. AP schemes have proven to be essential tools in capturing the correct asymptotic behavior of kinetic equations in stiff regimes (\cite{jin1999-AP,degond_asymptotic-preserving_2017,hu_chapter_2017}). Recent numerical studies have introduced APNNs~\cite{jin2023asymptotic}, a class of PINNs designed to maintain accuracy in the vanishing-Knudsen-number limit. Building on our earlier theoretical analysis, we now establish a rigorous generalization bound for APNNs. This result extends our convergence framework to a new class of kinetic models with multiscale structure.

% The asymptotic property of PINNs was proposed and numerically studied in APNNs (\cite{jin2023asymptotic}). As a generalization of the above analysis, we shall provide a rigorous proof of the asymptotic preserving property for APNNs. 

Specifically, we consider the following multiscale equation as in \cite{jin2023asymptotic},
\begin{equation}\la{eq:boltzmann_veps}
\veps\partial_t f^\veps + \xi\cdot\nabla_x f^\veps=\frac1\varepsilon Q(f^\veps,f^\veps), \quad\text{with } (t,x,\xi) \in \R^+\times\R^3\times\R^3,
\end{equation}
where the collision operator is linear and in the form of
\begin{equation}\label{eq:linearker}
 Q(f,f) = L(f) := \int_{\R^3} \alpha(\xi,\xi') \left( M(\xi)f(\xi') - M(\xi')f(\xi) \right) d\xi'.
\end{equation}
Here, $M(\xi)=(2\pi)^{3/2}\exp(-\abs{\xi}^2/2)$ is a global Maxwellian and $\alpha$ is a bounded function such that $0<\alpha_0\leq\alpha(\xi,\xi')\leq\alpha_1$ and $\alpha(\xi,\xi')=\alpha(\xi',\xi)$. 

The APNNs are designed based on the micro-macro decomposition, which we briefly summarize as below: One decomposes $f^\veps$ into the equilibrium part $\rho M$ and the non-equilibrium part $g$:
\begin{equation}\label{eq:decompostion}
 f^\veps(t,x,\xi) = \rho^\veps(t,x)M(\xi)+\veps g^\veps(t,x,\xi) ,
\end{equation}
with the constraint 
\begin{equation}\label{eq:decomposition-constraint}
 \rho^\veps(t,x)=\int_{\R^3} f^\veps(t,x,\xi) d\xi, \quad\text{or equivalently, }\int_{\R^3} g^\veps(t,x,\xi) d\xi \equiv 0.
\end{equation}
Then substituting \eqref{eq:decompostion} into \eqref{eq:boltzmann_veps} gives
\begin{equation}\label{eq:boltzmann-dcomp}
 \veps\pa_t\rho^\veps M + \veps^2 \pa_t g^\veps + \xi\cdot\nabla_x\rho^\veps M + \veps\xi\cdot\nabla_x g^\veps = L(g^\veps).
\end{equation}
Integrating over $\xi\in\R^3$ produces
\begin{equation}\label{eq:macro-part-bolt}
 \pa_t\rho^\veps + \nabla_x\cdot\int_{\R^3}\xi g^\veps d\xi = 0.
\end{equation}
Define the projection $\mathbf{P}_0: L^2(\R^3_\xi)\rightarrow\,\text{span}\{M\}$ w.r.t. the inner product $\displaystyle <f,g>=\int_{\R^3}\frac{fg}{M} d\xi$, and $\mathbf{P}_1=1-\mathbf{P}_0$, then from \eqref{eq:boltzmann-dcomp}, one obtains for $g^\veps$:
\begin{align}\label{eq:micro-part-bolt}
 \veps^2\pa_t g^\veps + \veps \mathbf{P}_1\left( \xi \cdot \nabla_x g^\veps \right) + \xi \cdot \nabla_x \rho^\veps M  = L(g^\veps).
\end{align}
The equations \eqref{eq:macro-part-bolt} and \eqref{eq:micro-part-bolt} together with the constraint \eqref{eq:decomposition-constraint} constitute the macro-micro decomposition formulation of \eqref{eq:boltzmann_veps}--\eqref{eq:linearker}.

As $\veps\rightarrow0$, formally, one obtains from \eqref{eq:micro-part-bolt},
\begin{equation}
     \xi \cdot \nabla_x \rho^0 M  = L(g^0),
\end{equation}
then by the solvability condition, 
\begin{equation}
    g^0(t,x,\xi) = h(\xi)\cdot\nabla_x\rho^0(t,x),
\end{equation}
where $h$ is the unique solution of $L(h)=\xi M(\xi)$ in $(\ker{L})^\perp$. Then the macro-component equation \eqref{eq:macro-part-bolt} converges to the diffusion equation 
\begin{equation}\label{eq:macro-diffusion}
    \pa_t\rho^0 + \nabla_x\cdot\int_{\R^3}\xi\otimes h(\xi) d\xi \cdot \nabla_x\rho^0(t,x) = 0.
\end{equation}

Figure~\ref{fig:APNN} is a reproduced plot from \cite{jin2023asymptotic} with slightly modified notations to be consistent with this paper. In this figure, $\mathcal{F}^\veps$ represents the solution to the multiscale equation \eqref{eq:boltzmann_veps}--\eqref{eq:linearker} (or equivalently \eqref{eq:macro-part-bolt}--\eqref{eq:micro-part-bolt}), and $\mathcal{F}^0$ represents the solution to the limit macroscopic model (diffusion equation) \eqref{eq:macro-diffusion}. Denote by $\mathcal{E}_G(\mathcal{F}^\veps)$ the generalization error of PINN associated with \eqref{eq:macro-part-bolt}--\eqref{eq:micro-part-bolt}, and then the asymptotic preserving property is to prove that when $\veps$ is small, the PINN solution obtained by minimizing $\mathcal{E}_G(\mathcal{F}^\veps)$ provides a good approximation of $\rho^0$, as long as $\mathcal{E}_G(\mathcal{F}^\veps)$ is small. Since it is shown that the distance between $\rho^\veps$ and $\rho^0$ is of $\mathcal{O}(\veps)$ (e.g., \cite{bardos_diffusion_nodate,poupaud_diffusion_1991,chai2018diffusion}), it suffices to prove that the total error of $\rho^\veps$ can be bounded by a multiple of $\mathcal{E}_G(\mathcal{F}^\veps)$.

\begin{figure}
    \centering
\begin{tikzpicture}
\node (1) at(0,0) {$\mathcal{F}^\varepsilon$};
\node (2) at(5,0) {$\mathcal{F}^0$};
\node (3) at(0,3) {$\mathcal{E}_G(\mathcal{F}^\varepsilon)$};
\node (4) at(5,3) {$\mathcal{E}_G(\mathcal{F}^0)$};
\draw[->] (1)--(2);
\draw[->] (3)--(4);
\draw[<-] (1)--(3);
\draw[<-] (2)--(4);
\node at(0.6,1.5) {$\mathcal{E}_G\rightarrow0$};
\node at(5.6,1.5) {$\mathcal{E}_G\rightarrow0$};
\node at(2.5,3.3) {$\varepsilon\rightarrow0$};
\node at(2.5,0.3) {$\varepsilon\rightarrow0$};
\end{tikzpicture}
    \caption{Interpretation of APNNs reproduced from \cite{jin2023asymptotic} with slightly modified notations to be consistent with this paper. Here, $\mathcal{F}^\veps$ represents the multiscale equation \eqref{eq:boltzmann_veps}--\eqref{eq:linearker}, and $\mathcal{F}^0$ represents the limit macroscopic model (diffusion equation) \eqref{eq:macro-diffusion}. Let $\mathcal{E}_G(\mathcal{F}^\veps)$ be the generalization error of PINN associated with \eqref{eq:macro-part-bolt}--\eqref{eq:micro-part-bolt}. The asymptotic preserving property aims to demonstrate that when $\veps$ is small, the PINN solution obtained by minimizing $\mathcal{E}_G(\mathcal{F}^\veps)$ can accurately approximate $\rho^0$, given that $\mathcal{E}_G(\mathcal{F}^\veps)$ is small. }
    \label{fig:APNN}
\end{figure}

\subsection{Residuals and Errors}

For a fixed small $\varepsilon > 0$ and a sufficiently large $R > 0$, we define the following residuals
\beg{align*}
\mathcal{R}_{1}[R; \varepsilon; \theta] (t, x) &= \pa_t \rho_{\theta}^{\varepsilon} + \na_x \cdot \int_{Q_R} \xi g_{\theta}^{\varepsilon} d\xi,
\\ \mathcal{R}_2 [R; \varepsilon;\theta] (t, x,\xi) &= \varepsilon^2 \pa_t g_{\theta}^{\varepsilon} 
+ \varepsilon \mathbf{P}_1 (\chi_{Q_R} \xi \cdot \na_x g_{\theta}^{\varepsilon}) + \xi \cdot \na_x \rho_{\theta}^{\varepsilon} M -  L(g_{\theta}^{\varepsilon} \chi_{Q_R}),
\\ \mathcal{R}_3 [\varepsilon;\theta](x) &= \rho_{\theta}^{\varepsilon} (0, x) - \rho^{\varepsilon}(0, x),
\\ \mathcal{R}_4 [\varepsilon; \theta] (x, \xi) &= g_{\theta}^{\varepsilon} (0, x, \xi) - g^{\varepsilon} (0, x, \xi),
\\ \mathcal{R}_5[R;\varepsilon;\theta] (t,x) &= \left(\int_{Q_R}g_{\theta}^{\varepsilon}  d\xi\right)^2 + \left|\na_x \int_{Q_R} g_{\theta}^{\varepsilon} d\xi\right|^2,
\end{align*} and we consider the following generalization error $\mathcal{E}_G[R; \varepsilon; \theta]$ and total error $\mathcal{E}_T [R; \varepsilon; \theta]$: 
\begin{equation}\label{eq:ap_total}
\beg{aligned}
\mathcal{E}_G[R; \varepsilon; \theta]^2 &= \int_{0}^{T} \int_{\TT^3} \mathcal{R}_1 [R; \varepsilon; \theta]^2 dxdt + \int_{0}^{T} \int_{\TT^3} \int_{Q_R} \frac{\mathcal{R}_2[R; \varepsilon; \theta]^2}{M} d\xi dx dt
\\&\quad\quad+ \int_{\TT^3} \mathcal{R}_3 [\varepsilon;\theta]^2 dx + \int_{\TT^3} \int_{Q_R} \mathcal{R}_4[\varepsilon; \theta]^2  d\xi dx + \int_{0}^{T} \int_{\TT^3} \mathcal{R}_5[R;\varepsilon;\theta]^2 dxdt ,
\\ \mathcal{E}_T [R; \varepsilon; \theta]^2 &= \int_{0}^{T} \int_{\TT^3} |\rho^{\varepsilon} - \rho_{\theta}^{\varepsilon}|^2 dxdt + \varepsilon^2 \int_{0}^{T} \int_{\TT^3} \int_{Q_R} \frac{|g^{\varepsilon} - g_{\theta}^{\varepsilon}|^2}{M} d\xi dxdt.
\end{aligned}    
\end{equation}
 The incorporation of the characteristic functions $\chi_{Q_R}$ and restriction to large cubes $Q_R$ are needed as 
$\xi \cdot \na_x g_{\theta}^{\varepsilon}$ and $g_{\theta}^{\varepsilon}$ do not necessarily belong to the domains of the operators $\mathbf{P}_1$ and $L$. 

\subsection{Preliminary Lemmas.} The following lemmas will be used to obtain good control of the total error $\mathcal{E}_T[R; \varepsilon;\theta]$ and rigorously prove the asymptotic preserving property. The incorporation of the charactersitic function $\chi_{Q_R}$ in the definition of the residuals gives rise to the need for deriving new estimates of the terms involving $\mathcal{P}_1$ and $L$.

\beg{lem} Let $R > 0$. For 
$\xi \in \R^3$, let 
\be 
\eta(\xi) = \int_{\RR^3} \alpha(\xi, \xi') M(\xi') d\xi'.
\ee Suppose there is a constant $c$ such that  
\be \label{sizea}
\frac{\|\alpha - c\|_{L^{\infty}(\RR^3 \times \RR^3)}}{\alpha_0} = {K} < 1.
\ee 
If $f \in L^2(\TT^3 \times Q_R)$, then it holds that 
\be \la{LL1}
\int_{\TT^3} \int_{Q_R} L(f\chi_{Q_R}) \frac{f}{M} d\xi dx   + (1-{K})  \|\sqrt{\frac{\eta}{M}} f\|_{L^2(\TT^3 \times Q_R)}^2 \le c \int_{\TT^3} \left(\int_{Q_R} f(\xi) d\xi \right)^2 dx. 
\ee  If $\tilde{f} \in L^2(\TT^3; L^1(\RR^3 \setminus Q_R))$, then it holds that 
\be \la{LL2}
\int_{\TT^3} \int_{Q_R} L(\tilde{f} \chi_{\R^3 \setminus Q_R}) \frac{f}{M} d\xi dx \le \frac{\alpha_1}{\sqrt{\alpha_0}} \|\tilde{f}\|_{L^2(\TT^3; L^1(\RR^3 \setminus Q_R))} \|\sqrt{\frac{\eta}{M}} f\|_{L^2(\TT^3 \times Q_R)} .
\ee 
\end{lem}

\beg{proof}.
We have 
\bes 
\beg{aligned}
&\int_{\TT^3} \int_{Q_R} L(f\chi_{Q_R})\frac{f}{M} d\xi dx
\\&\quad=\int_{\TT^3} \int_{Q_R} \int_{\RR^3} \alpha (\xi, \xi') \left(M(\xi) f(\xi')\chi_{Q_R}(\xi') - M(\xi') f(\xi)\chi_{Q_R}(\xi) \right) \frac{f(\xi)}{M(\xi)} d\xi'd\xi dx
\\&\quad= \int_{\TT^3} \int_{Q_R} \int_{Q_R} \alpha(\xi, \xi')  f(\xi)f(\xi') d\xi' d\xi dx - \int_{\TT^3} \int_{Q_R}\int_{\RR^3} \alpha(\xi, \xi') M(\xi') \frac{f(\xi)^2}{M(\xi)} d\xi' d\xi dx
\\&\quad= \int_{\TT^3} \int_{Q_R} \int_{Q_R} (\alpha (\xi, \xi') - c) f(\xi) f(\xi') d\xi' d\xi dx
+ c \int_{\TT^3} \int_{Q_R} \int_{Q_R} f(\xi) f(\xi') d\xi' d\xi dx 
\\&\quad\quad\quad\quad\quad\quad\quad\quad
- \left\|\sqrt{\frac{\eta}{M}} f\right\|_{L^2(\TT^3 \times Q_R)}^2
\\&\quad= \int_{\TT^3} \int_{Q_R} \int_{Q_R} (\alpha(\xi, \xi') - c) \sqrt{\frac{M(\xi)}{\eta(\xi)}} \sqrt{\frac{M(\xi')}{\eta(\xi')}} \sqrt{\frac{\eta(\xi)}{M(\xi)}}f(\xi)\sqrt{\frac{\eta(\xi')}{M(\xi')}}f(\xi') d\xi' d\xi dx  
\\&\quad\quad\quad\quad\quad\quad\quad\quad+c\int_{\TT^3} \left(\int_{Q_R} f(\xi) d\xi \right)^2 dx- \left\|\sqrt{\frac{\eta}{M}} f\right\|_{L^2(\TT^3 \times Q_R)}^2
\\&\quad\quad\le 
\|\alpha - c\|_{L^{\infty}(Q_R \times Q_R)} \left\| \sqrt{\frac{M}{\eta}}\right\|_{L^2(Q_R)}^2 \left\|\sqrt{\frac{\eta}{M}} f \right\|_{L^2(\TT^3 \times Q_R)}^2 
\\&\quad\quad\quad\quad\quad\quad\quad\quad+c\int_{\TT^3} \left(\int_{Q_R} f(\xi) d\xi \right)^2 dx
 - \left\|\sqrt{\frac{\eta}{M}} f\right\|_{L^2(\TT^3 \times Q_R)}^2,
\end{aligned}
\ees in view of the Cauchy-Schwarz inequality. 
Since $\eta(\xi) \ge \alpha_0$ for all $\xi \in \RR^3$, we estimate
\be 
\int_{Q_R} \frac{M}{\eta} d\xi \le \frac{1}{\alpha_0} \int_{\RR^3} M(\xi) d\xi = \frac{1}{\alpha_0},
\ee which yields
\be 
\int_{\TT^3} \int_{Q_R} L(f\chi_{Q_R}) \frac{f}{M} d\xi dx \le c
\int_{\TT^3} \left(\int_{Q_R} f(\xi) d\xi \right)^2 dx
 - (1-{K}) \left\|\sqrt{\frac{\eta}{M}} f\right\|_{L^2(\TT^3 \times Q_R)}^2,
\ee after making use of the size assumption \eqref{sizea}. This gives \eqref{LL1}. As for \eqref{LL2}, we have 
\bes 
\beg{aligned}
&\int_{\TT^3} \int_{Q_R} L(\tilde{f} \chi_{\R^3 \setminus Q_R}) \frac{f}{M} d\xi dx
\le \left\|\sqrt{\frac{\eta}{M}} f\right\|_{L^2(\TT^3 \times Q_R)} \left(\int_{\TT^3} \int_{Q_R} \frac{1}{\eta(\xi) M(\xi)} L(\tilde{f} \chi_{\RR^3 \setminus Q_R})^2 d\xi dx \right)^{\fr{1}{2}},
\end{aligned}
\ees by the Cauchy-Schwarz inequality. Since $\alpha \le \alpha_1$ and $\eta \ge \alpha_0$, we can bound 
\bes 
\beg{aligned}
&\int_{\TT^3} \int_{Q_R} \frac{1}{\eta(\xi) M(\xi)} L(\tilde{f} \chi_{\RR^3 \setminus Q_R})^2 d\xi dx 
\\&\quad\quad= \int_{\TT^3} \int_{Q_R} \frac{1}{\eta(\xi)M(\xi)}\left(\int_{\RR^3} \alpha (\xi, \xi') \left(M(\xi) \tilde{f}(\xi') \chi_{\RR^3 \setminus Q_R}(\xi') - M(\xi')\tilde{f}(\xi)\chi_{\RR^3 \setminus Q_R} (\xi) \right) d\xi' \right)^2 d\xi dx % cancellation used to obtain the following equality
\\&\quad\quad=  \int_{\TT^3} \int_{Q_R} \frac{1}{\eta (\xi) M(\xi)} \left(\int_{\RR^3 \setminus Q_R} \alpha (\xi, \xi') M(\xi) \tilde{f}(\xi')     d\xi' \right)^2 d\xi dx 
\\&\quad\quad\le \frac{\alpha_1^2}{\alpha_0}\int_{\TT^3} \int_{Q_R} M(\xi)  \|\tilde{f}\|_{L^1(\RR^3 \setminus Q_R)}^2  d\xi dx
%\\&\quad\quad\le \frac{\alpha_1^2}{\alpha_0} \|M\|_{L^1(Q_R)}\int_{\TT^3} \|\tilde{f}\|_{L^1(\RR^3 \setminus Q_R)}^2 dx
\le \frac{\alpha_1^2}{\alpha_0} \|M\|_{L^1(Q_R)} \|\tilde{f}\|_{L^2(\TT^3; L^1(\RR^3 \setminus Q_R))}^2.
\end{aligned}
\ees 
\end{proof}

\beg{rem}
As $\|\alpha - c\|_{L^{\infty}(\RR^3 \times \RR^3)} \le \max \left\{|\alpha_0 - c|, |\alpha_1 - c| \right\}$ for any constant $c$, the assumption \eqref{sizea} holds whenever there is a constant $c$ for which  $\max \left\{|\alpha_0 - c|, |\alpha_1 - c| \right\} < \alpha_0$. For instance, if $\alpha_0 = 1$ and $\alpha_1 = 1.5$, then one can choose $c = 0.75$ and obtain \eqref{sizea}. More generally, if $\alpha_1 < 2 \alpha_0$, then one can choose $c= \alpha_0$ and deduce \eqref{sizea}.  
\end{rem}

\beg{lem} Let $R>0$ and $f$ be a smooth function. Then it holds that
\be \label{cancel}
\beg{aligned}
\left|\int_{\TT^3} \int_{Q_R} \mathbf{P}_1 (\chi_{Q_R} \xi \cdot \na_x f) \frac{f}{M} d\xi dx \right| \le \frac{1- K}{4} \left\|\sqrt{\frac{\eta}{M}} f \right\|_{L^2(\TT^3 \times Q_R)}^2
+  C\int_{\TT^3} \left|\na_x \int_{Q_R} f d\xi \right|^2 dx,
\end{aligned}
\ee where $K$ is the constant defined in \eqref{sizea}. If $\tilde{f}$ is another smooth function, then we have
\be \la{cancel2}
\left|\int_{\TT^3} \int_{Q_R} \mathbf{P}_1 (\chi_{\RR^3 \setminus Q_R} \xi \cdot \na_x \tilde{f}) \frac{f}{M} d\xi dx \right|
\le  \frac{1- K}{4} \left\|\sqrt{\frac{\eta}{M}} f \right\|_{L^2(\TT^3 \times Q_R)}^2 + C\|\xi \cdot  \na_x \tilde{f}\|_{L^2(\TT^3; L^1(\RR^3 \setminus Q_R))}^2 .
\ee 
\end{lem}

\beg{proof}. On the one hand,
\be  
\beg{aligned}
\int_{\TT^3} \int_{Q_R} (\chi_{Q_R}(\xi) \xi \cdot \na_x f) \frac{f}{M} d\xi dx
= \int_{Q_R} \frac{1}{M} \int_{\TT^3} (\xi \cdot \na_x f)f dx d\xi = 0,
\end{aligned}
\ee because the inner spatial integral vanishes. On the other hand, 
\be 
\beg{aligned}
\left|\int_{\TT^3} \int_{Q_R} \mathbf{P}_0 (\chi_{Q_R} \xi \cdot \na_x f) \frac{{f}}{M} d\xi dx \right|
&= \left|\int_{\TT^3} \int_{Q_R} \left(\int_{\RR^3} \chi_{Q_R} \xi' \cdot \na_x f d\xi' \right) M \frac{{f}}{M} d\xi dx \right|
\\&= \left| \int_{\TT^3} \int_{Q_R} \left(\int_{Q_R} \xi' \cdot \na_x f d\xi' \right) {f} d\xi dx \right|
\\&= \left| \int_{\TT^3} \left(\int_{Q_R} \xi' \cdot \na_x f d\xi' \right) \left(\int_{Q_R} f d\xi \right) dx \right|
\\&= \left|-\int_{\TT^3} \left(\int_{Q_R} \xi' f d\xi'\right) \cdot \na_x \left(\int_{Q_R} f d\xi \right) dx \right|
\\&\le \left(\int_{\TT^3} \left(\int_{Q_R} \xi' f d\xi'\right)^2 dx \right)^{\fr{1}{2}} \left(\int_{\TT^3} \left|\na_x \int_{Q_R} f d\xi \right|^2 dx \right)^{\frac{1}{2}}.
\end{aligned}
\ee  Applying the Cauchy-Schwarz inequality in the $\xi$ variable gives
%\bes 
\beg{multline*}
\int_{\TT^3} \left(\int_{Q_R} \xi' f d\xi' \right)^2 dx
= \int_{\TT^3} \left(\int_{Q_R} \frac{\xi' \sqrt{M}}{\sqrt{\eta}} \sqrt{\frac{\eta}{M}} f d\xi' \right)^2 dx
\\\le \int_{\TT^3} \left\| \frac{\xi' \sqrt{M}}{\sqrt{\eta}} \right\|_{L^2(Q_R)}^2 \left\|\sqrt{\frac{\eta}{M}} f \right\|_{L^2(Q_R)}^2 dx
= \left\| \frac{\xi' \sqrt{M}}{\sqrt{\eta}} \right\|_{L^2(Q_R)}^2 \left\|\sqrt{\frac{\eta}{M}} f \right\|_{L^2(\TT^3 \times Q_R)}^2.
\end{multline*}
%\ees 
This gives rise to 
\bes 
\left|\int_{\TT^3} \int_{Q_R} \mathbf{P}_0 (\chi_{Q_R} \xi \cdot \na_x f) \frac{{f}}{M} d\xi dx \right|
\le \left\| \frac{\xi' \sqrt{M}}{\sqrt{\eta}} \right\|_{L^2(Q_R)} \left\|\sqrt{\frac{\eta}{M}} f \right\|_{L^2(\TT^3 \times Q_R)}\left(\int_{\TT^3} \left|\na_x \int_{Q_R} f d\xi \right|^2 \right)^{\frac{1}{2}}.
\ees An application of Young's inequality for products yields
\bes 
\beg{aligned}
\left|\int_{\TT^3} \int_{Q_R} \mathbf{P}_0 (\chi_{Q_R} \xi \cdot \na_x f) \frac{{f}}{M} d\xi dx \right|
&\le \frac{1- K}{4} \left\|\sqrt{\frac{\eta}{M}} f \right\|_{L^2(\TT^3 \times Q_R)}^2
\\&\quad+ C_K \left\| \frac{\xi' \sqrt{M}}{\sqrt{\eta}} \right\|_{L^2(Q_R)}^2 \int_{\TT^3} \left|\na_x \int_{Q_R} f d\xi \right|^2 dx. 
\end{aligned}
\ees Since 
\bes 
\int_{Q_R} \frac{\xi'^2 M(\xi')}{\eta} d\xi' \le \frac{1}{\alpha_0} \int_{\RR^3} \xi'^2 M(\xi') d\xi' \le C(\alpha_0),
\ees we obtain the desired estimate \eqref{cancel}. As for \eqref{cancel2}, the cancellation law
\bes 
\beg{aligned}
\int_{\TT^3} \int_{Q_R} \chi_{\RR^3 \setminus Q_R} (\xi \cdot \na_x \tilde{f}) \frac{f}{M} d\xi dx  = \int_{\TT^3} \int_{\R^3} \chi_{Q_R} \chi_{\RR^3 \setminus Q_R} (\xi \cdot \na_x \tilde{f}) \frac{f}{M} d\xi dx = 0
\end{aligned}
\ees holds due to the vanishing of the function $\chi_{\RR^3 \setminus Q_R} \chi_{Q_R} = 0$. Moreover,
\be 
\beg{aligned}
&\left|\int_{\TT^3} \int_{Q_R} \mathbf{P}_0(\chi_{\RR^3 \setminus Q_R} \xi \cdot \na_x \tilde{f}) \frac{f}{M} d\xi dx \right|
= \left|\int_{\TT^3} \left(\int_{\RR^3 \setminus Q_R} \xi' \cdot \na_x \tilde{f} d\xi' \right)\left(\int_{Q_R} f d\xi \right) dx \right|
\\&\quad\quad\le \int_{\TT^3}
\left\| \xi \na_x \tilde{f} \right\|_{L^1(\RR^3 \setminus Q_R)}  \left\| \frac{\sqrt{\eta}}{\sqrt{M}}f \right\|_{L^2(Q_R)} \left\|\frac{\sqrt{M}}{\sqrt{\eta}} \right\|_{L^2(Q_R)} dx
\\&\quad\quad\le \left\|\frac{\sqrt{M}}{\sqrt{\eta}} \right\|_{L^2(\RR^3)}  \left\| \frac{\sqrt{\eta}}{\sqrt{M}}f \right\|_{L^2(\TT^3 \times Q_R)} \left\| \xi \na_x \tilde{f} \right\|_{L^2(\TT^3; L^1(\RR^3 \setminus Q_R))},
\end{aligned}
\ee yielding \eqref{cancel2} after making use of Young's inequality. 
\end{proof}

\subsection{Asymptotic Preserving Property.} In this subsection, we prove the asymptotic preserving property.

\beg{Thm} Fix $\varepsilon > 0$. Let $T>0$ be an arbitrary positive time. Let $(\rho^{\varepsilon}, g^{\varepsilon})$ be a solution to \eqref{eq:macro-part-bolt}--\eqref{eq:micro-part-bolt}. Suppose there is a radius $R_0 >0$ and a positive constant $\delta(\epsilon)$  with $\delta(\epsilon) \rightarrow 0$ as $\epsilon \rightarrow 0$ such that 
\be 
\sup\limits_{R \ge R_0} \int_{0}^{T} \left\{\|g^{\varepsilon}\|_{L^2(\TT^3; L^1(\RR^3 \setminus Q_R))}^2 +  \|\xi \cdot \na_x g^{\varepsilon}\|_{L^2(\TT^3; L^1(\RR^3 \setminus Q_R))}^2 +  \varepsilon  \|\na_x g^{\varepsilon}\|_{L^2(\TT^3; L^1(\RR^3 \setminus Q_R))}^2 \right\} dt\le \delta(\epsilon).
\ee 
Let $R \ge R_0$. Let $(\rho_{\theta}^{\varepsilon}, g_{\theta}^{\varepsilon})$ be \texttt{tanh} neural networks approximating $(\rho^{\varepsilon}, g^{\varepsilon})$. Under assumption \eqref{sizea}, there exists a positive constant $C$ depending only on $T$ such that 
\be \label{APPst}
\mathcal{E}_T [R; \varepsilon; \theta]^2
\le C\left(\mathcal{E}_G [R; \varepsilon; \theta]^2 +\delta(\epsilon) \right).
\ee 
  \end{Thm}

Remark that the way the total error is defined in \eqref{eq:ap_total}, along with the estimate \eqref{APPst}, shows that the macroscopic part $\rho^\varepsilon_\theta$ approximates $\rho^\varepsilon$ well when the generalization error is small. However, resolving $\varepsilon$ is necessary to guarantee the convergence of $g^{\varepsilon}$. This is precisely the asymptotic preserving property.
 
  \begin{proof}.
The differences $\rho^{\varepsilon} - \rho_{\theta}^{\varepsilon}$ and $g^{\varepsilon} - g_{\theta}^{\varepsilon}$ evolve according to 
\be \label{appmain1}
\pa_t (\rho^{\varepsilon} - \rho_{\theta}^{\varepsilon})
+ \na_x \cdot \int_{Q_R} \xi \cdot (g^{\varepsilon} - g_{\theta}^{\varepsilon}) d\xi 
+ \na_x \cdot \int_{\RR^3 \setminus Q_R} \xi \cdot g^{\varepsilon} d\xi = - \mathcal{R}_1,
\ee and
\be \label{appmain2}
\beg{aligned}
&\varepsilon^2 \pa_t (g^{\varepsilon} - g_{\theta}^{\varepsilon}
) + \varepsilon \mathbf{P}_1 \left(\chi_{Q_R} \xi \cdot \na_x (g^{\varepsilon} - g_{\theta}^{\varepsilon}) \right)
+ \varepsilon \mathbf{P}_1 (\chi_{\RR^3 \setminus Q_R} \xi \cdot \na_x g^{\varepsilon}) + \xi \cdot \na_x(\rho^{\varepsilon} - \rho_{\theta}^{\varepsilon}) M 
\\&\quad\quad\quad\quad- L((g^{\varepsilon} - g_{\theta}^{\varepsilon})\chi_{Q_R}) - L(g^{\varepsilon} \chi_{\RR^3 \setminus Q_R}) = - \mathcal{R}_2,
\end{aligned}
\ee respectively. We multiply \eqref{appmain1} by $\rho^{\varepsilon} - \rho_{\theta}^{\varepsilon}$ and integrate over $\TT^3$. We multiply \eqref{appmain2} by $\frac{g^{\varepsilon} - g_{\theta}^{\varepsilon}}{M}$ and integrate over $\TT^3 \times Q_R$. Then we add the resulting energy equalities. In view of the cancellation law
\be 
\int_{\TT^3} \na_x \cdot \int_{Q_R} \xi \cdot (g^{\varepsilon} - g_{\theta}^{\varepsilon}) (\rho^{\varepsilon} - \rho_{\theta}^{\varepsilon}) d\xi dx + \int_{\TT^3} \int_{Q_R} \xi \cdot \na_x (\rho^{\varepsilon} - \rho_{\theta}^{\varepsilon})M \frac{g^{\varepsilon} - g_{\theta}^{\varepsilon}}{M} d\xi dx = 0,
\ee we obtain 
\be 
\beg{aligned}
&\frac{1}{2} \frac{d}{dt} \left(\|\rho^{\varepsilon} - \rho_{\theta}^{\varepsilon} \|_{L^2(\TT^3)}^2 + {\varepsilon^2} \left\|\frac{g^{\varepsilon} - g_{\theta}^{\varepsilon}}{\sqrt{M}} \right\|_{L^2(\TT^3 \times Q_R)}^2 \right)
= -\int_{\TT^3} \int_{\RR^3 \setminus Q_R} \xi \cdot \na_x g^{\varepsilon} (\rho^{\varepsilon} - \rho_{\theta}^{\varepsilon}) d\xi dx 
\\&\quad- \int_{\TT^3} \mathcal{R}_1  (\rho^{\varepsilon} - \rho_{\theta}^{\varepsilon}) dx
-\varepsilon \int_{\TT^3} \int_{Q_R} \mathbf{P}_1 (\chi_{\RR^3 \setminus Q_R} \xi \cdot \na_x (g^{\varepsilon} - g_{\theta}^{\varepsilon})) \frac{g^{\varepsilon} - g_{\theta}^{\varepsilon}}{M} d\xi dx 
\\&\quad\quad- \varepsilon \int_{\TT^3} \int_{Q_R} \mathbf{P}_1 (\chi_{\RR^3 \setminus Q_R} \xi \cdot \na_x g^{\varepsilon}) \frac{g^{\varepsilon} - g_{\theta}^{\varepsilon}}{M} d\xi dx 
+\int_{\TT^3} \int_{Q_R} L((g^{\varepsilon} - g_{\theta}^{\varepsilon})\chi_{Q_R}) \frac{g^{\varepsilon} - g_{\theta}^{\varepsilon}}{M} d\xi dx 
\\&\quad\quad\quad\quad+ \int_{\TT^3} \int_{Q_R} L(g^{\varepsilon}\chi_{\RR^3 \setminus Q_R}) \frac{g^{\varepsilon} - g_{\theta}^{\varepsilon}}{M} d\xi dx
- \int_{\TT^3} \int_{Q_R} \mathcal{R}_2 \frac{g^{\varepsilon} - g_{\theta}^{\varepsilon}}{M} d\xi dx.
\end{aligned}
\ee We bound
\be 
\left|\int_{\TT^3} \mathcal{R}_1  (\rho^{\varepsilon} - \rho_{\theta}^{\varepsilon}) dx \right| \le C\|\mathcal{R}_1\|_{L^2(\TT^3)}^2 + \|\rho^{\varepsilon} - \rho_{\theta}^{\varepsilon} \|_{L^2(\TT^3)}^2,
\ee and 
\be 
\left| \int_{\TT^3} \int_{Q_R} \mathcal{R}_2 \frac{g^{\varepsilon} - g_{\theta}^{\varepsilon}}{M} d\xi dx \right|
\le C\|\frac{\mathcal{R}_2}{\sqrt{\eta M}} \|_{L^2(\TT^3 \times Q_R)}^2 + \frac{1-K}{16} \|\sqrt{\frac{\eta}{M}} (g^{\varepsilon} - g_{\theta}^{\varepsilon})\|_{L^2(\TT^3 \times Q_R)}^2,
\ee using the Cauchy-Schwarz and Young inequalities. In view of \eqref{LL1}, \eqref{LL2}, \eqref{cancel}, and \eqref{cancel2}, we estimate
\bes
\beg{aligned}
&\int_{\TT^3} \int_{Q_R} L((g^{\varepsilon} - g_{\theta}^{\varepsilon})\chi_{Q_R}) \frac{g^{\varepsilon} - g_{\theta}^{\varepsilon}}{M} d\xi dx   + (1-K)  \|\sqrt{\frac{\eta}{M}} (g^{\varepsilon} - g_{\theta}^{\varepsilon})\|_{L^2(\TT^3 \times Q_R)}^2 \\ &\hspace{15em} \le C\int_{\TT^3} \left(\int_{Q_R} (g^{\varepsilon} - g_{\theta}^{\varepsilon}) d\xi \right)^2 dx, \\
&\int_{\TT^3} \int_{Q_R} L(g^{\varepsilon} \chi_{\R^3 \setminus Q_R}) \frac{g^{\varepsilon} - g_{\theta}^{\varepsilon}}{M} d\xi dx \le C\|g^{\varepsilon} \|_{L^2(\TT^3; L^1(\RR^3 \setminus Q_R))} \|\sqrt{\frac{\eta}{M}} (g^{\varepsilon} - g_{\theta}^{\varepsilon})\|_{L^2(\TT^3 \times Q_R)}, \\
&\left|\int_{\TT^3} \int_{Q_R} \mathbf{P}_1 (\chi_{Q_R} \xi \cdot \na_x (g^{\varepsilon} - g_{\theta}^{\varepsilon})) \frac{g^{\varepsilon} - g_{\theta}^{\varepsilon}}{M} d\xi dx \right| \\ &\hspace{15em} 
\le \frac{1- K}{4} \left\|\sqrt{\frac{\eta}{M}} (g^{\varepsilon} - g_{\theta}^{\varepsilon}) \right\|_{L^2(\TT^3 \times Q_R)}^2
+  C\int_{\TT^3} \left|\na_x \int_{Q_R} (g^{\varepsilon} - g_{\theta}^{\varepsilon}) d\xi \right|^2 dx,
\end{aligned}
\ees and
\bes 
\beg{aligned}
\left|\int_{\TT^3} \int_{Q_R} \mathbf{P}_1 (\chi_{\RR^3 \setminus Q_R} \xi \cdot \na_x g^{\varepsilon} ) \frac{g^{\varepsilon} - g_{\theta}^{\varepsilon}}{M} d\xi dx \right|
&\le  \frac{1- K}{4} \left\|\sqrt{\frac{\eta}{M}}  (g^{\varepsilon} - g_{\theta}^{\varepsilon})\right\|_{L^2(\TT^3 \times Q_R)}^2 
\\&\quad\quad+ C\|\xi \cdot  \na_x g^{\varepsilon} \|_{L^2(\TT^3; L^1(\RR^3 \setminus Q_R))}^2.
\end{aligned}
\ees Consequently, we obtain the differential inequality
\bes 
\beg{aligned}
&\frac{1}{2} \frac{d}{dt} \left(\|\rho^{\varepsilon} - \rho_{\theta}^{\varepsilon} \|_{L^2(\TT^3)}^2 + {\varepsilon^2} \left\|\frac{g^{\varepsilon} - g_{\theta}^{\varepsilon}}{\sqrt{M}} \right\|_{L^2(\TT^3 \times Q_R)}^2 \right)
+ \frac{1-K}{4}\left\|\sqrt{\frac{\eta}{M}}  (g^{\varepsilon} - g_{\theta}^{\varepsilon})\right\|_{L^2(\TT^3 \times Q_R)}^2 
\\&\quad\quad\le C\|\rho^{\varepsilon} - \rho_{\theta}^{\varepsilon} \|_{L^2(\TT^3)}^2 + C\left(\|\mathcal{R}_1\|_{L^2(\TT^3)}^2 + \|\frac{\mathcal{R}_2}{\sqrt{M}}\|_{L^2(\TT^3 \times Q_R)}^2 \right) 
\\&\quad\quad\quad\quad+C\left((1+\varepsilon) \|\xi \cdot  \na_x g^{\varepsilon} \|_{L^2(\TT^3; L^1(\RR^3 \setminus Q_R))}^2 + \|g^{\varepsilon} \|_{L^2(\TT^3; L^1(\RR^3 \setminus Q_R))} ^2\right)
\\&\quad\quad\quad\quad\quad\quad+ C\left(\int_{\TT^3} \left(\int_{Q_R} (g^{\varepsilon} - g_{\theta}^{\varepsilon}) d\xi \right)^2 dx +  \varepsilon \int_{\TT^3} \left|\na_x \int_{Q_R} (g^{\varepsilon} - g_{\theta}^{\varepsilon}) d\xi \right|^2 dx\right),
\end{aligned}
\ees 
after applying Young's inequality for products. Using the relation 
\bes 
\int_{Q_R} g^{\varepsilon} d\xi + \int_{\RR^3 \setminus Q_R} g^{\varepsilon} d\xi = 0, 
\ees we have \bes
\left(\int_{Q_R} (g^{\varepsilon} - g_{\theta}^{\varepsilon}) d\xi \right)^2 = \left(-\int_{\RR^3 \setminus Q_R} g^{\varepsilon} d\xi - \int_{Q_R} g_{\theta}^{\varepsilon} d\xi \right)^2 
\le C\left(\int_{\RR^3 \setminus Q_R} g^{\varepsilon} d\xi  \right)^2 + C\left(\int_{Q_R} g_{\theta}^{\varepsilon} d\xi \right)^2,
\ees and 
\bes
\left|\na_x \int_{Q_R} (g^{\varepsilon} - g_{\theta}^{\varepsilon}) d\xi \right|^2  
\le C\left|\int_{\RR^3 \setminus Q_R} \na_x g^{\varepsilon} d\xi  \right|^2 + C\left|\int_{Q_R} \na_x g_{\theta}^{\varepsilon} d\xi \right|^2.
\ees Therefore, we obtain 
\bes 
\beg{aligned}
& \frac{d}{dt} \left(\|\rho^{\varepsilon} - \rho_{\theta}^{\varepsilon} \|_{L^2(\TT^3)}^2 + {\varepsilon^2} \left\|\frac{g^{\varepsilon} - g_{\theta}^{\varepsilon}}{\sqrt{M}} \right\|_{L^2(\TT^3 \times Q_R)}^2 \right) - C\left(\|\rho^{\varepsilon} - \rho_{\theta}^{\varepsilon} \|_{L^2(\TT^3)}^2 + \varepsilon^2 \left\|\frac{g^{\varepsilon} - g_{\theta}^{\varepsilon}}{\sqrt{M}} \right\|_{L^2(\TT^3 \times Q_R)}^2 \right)
\\&\quad\quad\le  C\left(\|\mathcal{R}_1\|_{L^2(\TT^3)}^2 + \|\frac{\mathcal{R}_2}{\sqrt{M}}\|_{L^2(\TT^3 \times Q_R)}^2 + \|\mathcal{R}_5 \|_{L^1(\TT^3)} \right) 
\\&\quad\quad\quad\quad+C\left(\|\xi \cdot  \na_x g^{\varepsilon} \|_{L^2(\TT^3; L^1(\RR^3 \setminus Q_R))}^2 + \|g^{\varepsilon} \|_{L^2(\TT^3; L^1(\RR^3 \setminus Q_R))} ^2 +    \varepsilon \|\na_x g^{\varepsilon}\|_{L^2(\TT^3; L^1(\RR^3 \setminus Q_R))}^2\right).
\end{aligned}
\ees By Gronwall's inequality, we deduce that \eqref{APPst} holds. 
  \end{proof} 

%%%%
%NUMERICAL EXAMPLE
%%%%
\section{Numerical experiments}\label{sec:numerical}
In this section, we present a numerical illustration to support the theoretical results and demonstrate the feasibility of using PINNs for nonlocal kinetic models. For simplicity, we consider a linearized version of the Boltzmann equation near a global Maxwellian in one spatial dimension. While this example does not capture the full complexity of the nonlinear Boltzmann dynamics, it offers a controlled setting to test the implementation of the proposed mollified PINN framework and to validate the convergence behavior predicted by our analysis:
\begin{equation}
	\partial_t f + \xi \cdot \nabla_x f = 
	\int_{\R} \left( M(\xi)f(\xi') - M(\xi')f(\xi) \right) d\xi',
	\quad \text{with } (t, x, \xi) \in \R^+ \times \mathbb{T} \times \R,
\end{equation}
where the global Maxwellian is given by $M(\xi) = (2\pi)^{-1/2} \exp(-25\xi^2/2)$.
The initial condition is specified as
\begin{equation}
	f(0,x,\xi) = \left( 1 + \cos(4\pi x) \right) \exp(-50\xi^2).
\end{equation}

We use a fully-connected neural network with $\texttt{m}$ layers, each containing $\texttt{n}$ neurons (denoted by $\texttt{FCmXn}$). The \texttt{tanh} activation function is used for all layers. The training process utilizes the Adam optimizer, and the loss function is constructed based on the PDE residual, the initial residual, and the boundary residual to ensure consistency with the governing equations.

The reference solution is obtained using a high-resolution spectral method, which serves as the benchmark for evaluating the accuracy of the PINN solutions.

\begin{figure}[htp]
	\centering
	\includegraphics[width=\textwidth]{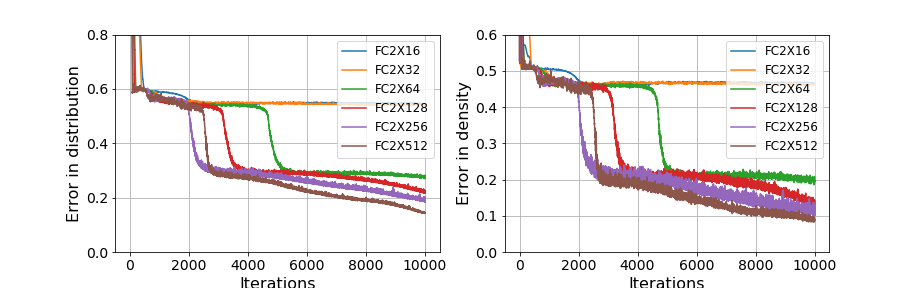}\\
	\includegraphics[width=\textwidth]{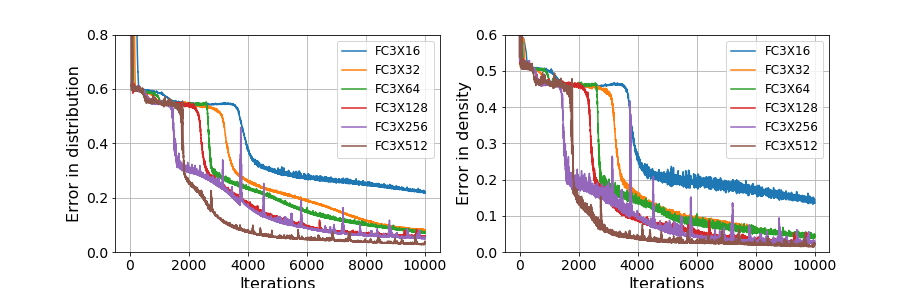}\\
	\includegraphics[width=\textwidth]{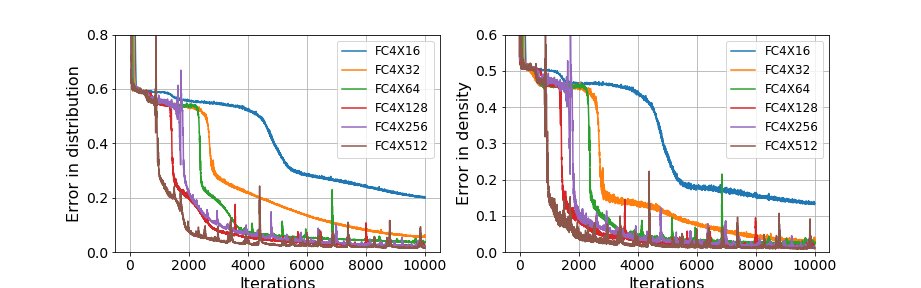}
	\caption{Comparison of the PINN solution and the reference solution for the Boltzmann equation. The left column displays the errors in the distribution functions, while the right column shows the errors in the density functions. From top to bottom, the results correspond to fully-connected neural networks with 2, 3, and 4 layers, each tested with layer widths of 16, 32, 64, 128, 256, and 512 neurons.}
	\label{fig:comparison}
\end{figure}

Figure~\ref{fig:comparison} illustrates the performance of the PINN solutions compared to the reference solutions for both the distribution function $f$ and the density function $\rho$ as the iterations progress. The errors decrease as the network depth and width increase, demonstrating the capability of PINNs to approximate the solution with high accuracy.

\section{Conclusion}\label{sec:conclusion}
In this paper, motivated by the successful application of physics-informed neural networks (PINNs) to solve Boltzmann-type equations \cite{jin2023asymptotic}, we have provided a rigorous error analysis for PINNs in approximating the solution of the Boltzmann equation near a global Maxwellian. The primary challenge in this analysis arises from the nonlocal quadratic interaction term, which is defined in the unbounded domain of velocity space. To address this, we incorporated a truncation function, necessitating the use of delicate analysis techniques. Furthermore, we extended our analysis to demonstrate the asymptotic preserving property of PINNs when employing micro-macro decomposition-based neural networks. This generalization not only reinforces the robustness of PINNs in handling Boltzmann-type equations but also underscores their potential for broader applications in kinetic theory and related fields. Our findings contribute to the theoretical foundations of PINNs, providing valuable insights into their accuracy and effectiveness in solving high-dimensional and complex PDEs. We hope this work may pave the way for future research to explore and optimize PINNs for a wider range of applications in computational physics and beyond.

\section*{Acknowledgements}
L.C. was partially supported by the National Key R\&D Program of China No. 2021YFA1003001 and the NSFC grant No. 12271537.  R.H. was partially supported by the ONR grant under \#N00014-24-1-2432, the Simons Foundation (MP-TSM-00002783) and the NSF grant DMS-2420988. X.Y. gratefully acknowledges partial support from the NSF grant DMS-2109116. Part of this work was done during L.C.'s visit to the Department of Mathematics at the University of California, Santa Barbara. L.C. would like to express gratitude to the department for their hospitality. The authors thank the anonymous referee for valuable suggestions, which have helped improve the clarity and quality of the paper.

%USE THE BELOW OPTIONS IN CASE YOU NEED AUTHOR YEAR FORMAT.
\bibliographystyle{plain}
\bibliography{reference}

@article{jin1999-AP,
author = {Jin, Shi},
title = {Efficient Asymptotic-Preserving ({AP}) Schemes For Some Multiscale Kinetic Equations},
journal = {SIAM Journal on Scientific Computing},
volume = {21},
number = {2},
pages = {441-454},
year = {1999}
}

@article{degond_asymptotic-preserving_2017,
	title = {Asymptotic-Preserving methods and multiscale models for plasma physics},
	volume = {336},
	issn = {0021-9991},
	journal = {Journal of Computational Physics},
	author = {Degond, Pierre and Deluzet, Fabrice},
	month = may,
	year = {2017},
	pages = {429--457},
}

@incollection{hu_chapter_2017,
	series = {Handbook of {Numerical} {Methods} for {Hyperbolic} {Problems}},
	title = {Chapter 5 - {Asymptotic}-{Preserving} {Schemes} for {Multiscale} {Hyperbolic} and {Kinetic} {Equations}},
	volume = {18},
	booktitle = {Handbook of {Numerical} {Analysis}},
	publisher = {Elsevier},
	author = {Hu, J. and Jin, S. and Li, Q.},
	editor = {Abgrall, Rémi and Shu, Chi-Wang},
	month = jan,
	year = {2017},
	pages = {103--129},
}

@article{xiao_using_2021,
	title = {Using neural networks to accelerate the solution of the {Boltzmann} equation},
	volume = {443},
	issn = {00219991},
	language = {en},
	journal = {Journal of Computational Physics},
	author = {Xiao, Tianbai and Frank, Martin},
	month = oct,
	year = {2021},
	pages = {110521},
}

@article{han_uniformly_2019,
	title = {Uniformly accurate machine learning-based hydrodynamic models for kinetic equations},
	volume = {116},
	number = {44},
	journal = {Proceedings of the National Academy of Sciences},
	author = {Han, Jiequn and Ma, Chao and Ma, Zheng and E, Weinan},
	month = oct,
	year = {2019},
	note = {Publisher: Proceedings of the National Academy of Sciences},
	pages = {21983--21991},
}

@article{li_physics-informed_2023,
	title = {Physics-Informed Deep Learning for Solving Coupled Electron and Phonon {Boltzmann} Transport Equations},
	volume = {19},
	issn = {2331-7019},
	language = {en},
	number = {6},
	journal = {Physical Review Applied},
	author = {Li, Ruiyang and Lee, Eungkyu and Luo, Tengfei},
	month = jun,
	year = {2023},
	pages = {064049},
}

@article{li_physics-informed_2022,
	title = {Physics-informed deep learning for solving phonon {Boltzmann} transport equation with large temperature non-equilibrium},
	volume = {8},
	issn = {2057-3960},
	language = {en},
	number = {1},
	journal = {npj Computational Materials},
	author = {Li, Ruiyang and Wang, Jian-Xun and Lee, Eungkyu and Luo, Tengfei},
	month = feb,
	year = {2022},
	pages = {29},
}

@article{lou_physics-informed_2021,
	title = {Physics-informed neural networks for solving forward and inverse flow problems via the {Boltzmann}-{BGK} formulation},
	volume = {447},
	issn = {00219991},
	language = {en},
	journal = {Journal of Computational Physics},
	author = {Lou, Qin and Meng, Xuhui and Karniadakis, George Em},
	month = dec,
	year = {2021},
	pages = {110676},
}

@article{bardos_diffusion_nodate,
	title = {DIFFUSION APPROXIMATION AND COMPUTATION OF THE CRITICALSIZE},
	language = {en},
	author = {Bardos, C and Sentis, R Santosand R},
    journal = {Trans. Amer. Math. Soc. },
    volume = {284},
    year = {1984},
    pages = {617--649},
}

@article{poupaud_diffusion_1991,
	title = {Diffusion approximation of the linear semiconductor {Boltzmann} equation: analysis of boundary layers},
	volume = {4},
	number = {4},
	journal = {Asymptotic Analysis},
	author = {Poupaud, F.},
	year = {1991},
	pages = {293--317},
}

@article{chai2018diffusion,
  title={{Diffusion limit of the Boltzmann--Landau--Lifshitz--Gilbert system in ferromagnetic materials}},
  author={Chai, Lihui and Garc{\'\i}a-Cervera, Carlos J and Yang, Xu},
  journal={Communications in Mathematical Sciences},
  volume={16},
  number={4},
  pages={1157--1167},
  year={2018},
  publisher={International Press of Boston}
}

@article{de2021approximation,
	title={On the approximation of functions by tanh neural networks},
	author={De Ryck, Tim and Lanthaler, Samuel and Mishra, Siddhartha},
	journal={Neural Networks},
	volume={143},
	pages={732--750},
	year={2021}
}

@article{de2022generic,
  title={Generic bounds on the approximation error for physics-informed (and) operator learning},
  author={De Ryck, Tim and Mishra, Siddhartha},
  journal={Advances in Neural Information Processing Systems},
  volume={35},
  pages={10945--10958},
  year={2022}
}

@article{de2022error,
  title={Error analysis for physics-informed neural networks ({PINN}s) approximating {K}olmogorov {PDE}s},
  author={De Ryck, Tim and Mishra, Siddhartha},
  journal={Advances in Computational Mathematics},
  volume={48},
  number={6},
  pages={79},
  year={2022},
  publisher={Springer}
}

@article{de2024error,
	title={Error estimates for physics-informed neural networks approximating the {N}avier--{S}tokes equations},
	author={De Ryck, Tim and Jagtap, Ameya D and Mishra, Siddhartha},
	journal={IMA Journal of Numerical Analysis},
	volume={44},
	number={1},
	pages={83--119},
	year={2024}
}

@article{mishra2022estimates,
  title={Estimates on the generalization error of physics-informed neural networks for approximating a class of inverse problems for {PDEs}},
  author={Mishra, Siddhartha and Molinaro, Roberto},
  journal={IMA Journal of Numerical Analysis},
  volume={42},
  number={2},
  pages={981--1022},
  year={2022}
}

@article{mishra2023estimates,
  title={Estimates on the generalization error of physics-informed neural networks for approximating {PDEs}},
  author={Mishra, Siddhartha and Molinaro, Roberto},
  journal={IMA Journal of Numerical Analysis},
  volume={43},
  number={1},
  pages={1--43},
  year={2023}
}

@article{jiao2021rate,
  title={A rate of convergence of Physics Informed Neural Networks for the linear second order elliptic {PDEs}},
  author={Jiao, Yuling and Lai, Yanming and Li, Dingwei and Lu, Xiliang and Wang, Fengru and Wang, Yang and Yang, Jerry Zhijian},
  journal={Communications in Computational Physics},
  volume={31},
  number={4},
  pages={1272--1295},
  year={2022}
}

@Article{jiao2023rate,
  author = {Jiao, Yuling and Yang, Jerry Zhijian and Yuan, Cheng and Zhou, Junyu},
  title = {A Rate of Convergence of Weak Adversarial Neural Networks for the Second Order Parabolic PDEs},
  journal = {Communications in Computational Physics},
  year = {2023},
  volume = {34},
  number = {3},
  pages = {813--836}
}

@Article{jiao2023improved,
  author = {Jiao, Yuling and Lu, Xiliang and Yang, Jerry Zhijian and Yuan, Cheng and Zhang, Pingwen},
  title = {Improved Analysis of {PINNs}: Alleviate the {CoD} for Compositional Solutions},
  journal = {Annals of Applied Mathematics},
  year = {2023},
  volume = {39},
  number = {3},
  pages = {239--263}
}

@article{guo2022mc,
  title = {{Monte Carlo fPINNs}: Deep learning method for forward and inverse problems involving high dimensional fractional partial differential equations},
  journal = {Computer Methods in Applied Mechanics and Engineering},
  volume = {400},
  pages = {115523},
  year = {2022},
  author = {Ling Guo and Hao Wu and Xiaochen Yu and Tao Zhou}
}

@article{pang2020npinn,
  title = {{nPINNs}: Nonlocal physics-informed neural networks for a parametrized nonlocal universal {Laplacian} operator. {Algorithms and applications}},
  journal = {Journal of Computational Physics},
  volume = {422},
  pages = {109760},
  year = {2020},
  author = {G. Pang and M. D'Elia and M. Parks and G.E. Karniadakis}
}

@article{biswas2022error,
	title={Error estimates for deep learning methods in fluid dynamics},
	author={Biswas, Animikh and Tian, Jing and Ulusoy, Suleyman},
	journal={Numerische Mathematik},
	volume={151},
	number={3},
	pages={753--777},
	year={2022},
	publisher={Springer}
}

@article{hu2023higher,
	title={Higher-order error estimates for physics-informed neural networks approximating the primitive equations},
	author={Hu, Ruimeng and Lin, Quyuan and Raydan, Alan and Tang, Sui},
	journal={Partial Differential Equations and Applications},
	volume={4},
	number={4},
	pages={34},
	year={2023}
}

@article{abdo2024accuracy,
	title={Accuracy Analysis of Physics-Informed Neural Networks for Approximating the Critical {SQG} Equation},
	author={Abdo, Elie and Hu, Ruimeng and Lin, Quyuan},
	journal={arXiv preprint arXiv:2401.10879},
	year={2024}
}

@article{jin2023asymptotic,
  title={Asymptotic-preserving neural networks for multiscale time-dependent linear transport equations},
  author={Jin, Shi and Ma, Zheng and Wu, Keke},
  journal={Journal of Scientific Computing},
  volume={94},
  number={3},
  pages={57},
  year={2023}
}

@article{poette2019gpc,
  title={{A gPC-intrusive Monte-Carlo scheme for the resolution of the uncertain linear Boltzmann equation}},
  author={Po{\"e}tte, Ga{\"e}l},
  journal={Journal of Computational Physics},
  volume={385},
  pages={135--162},
  year={2019},
  publisher={Elsevier}
}

@article{poette2022numerical,
  title={{Numerical analysis of the Monte-Carlo noise for the resolution of the deterministic and uncertain linear Boltzmann equation (comparison of non-intrusive gPC and MC-gPC)}},
  author={Po{\"e}tte, Ga{\"e}l},
  journal={Journal of Computational and Theoretical Transport},
  pages={1--53},
  year={2022},
  publisher={Taylor \& Francis}
}

@incollection{Cerci,
  title={The {B}oltzmann equation},
  author={Cercignani, Carlo},
  booktitle={The Boltzmann equation and its applications},
  pages={40--103},
  year={1988},
  publisher={Springer}
}

@book{Jin-Pareschi-Book,
  title={Uncertainty quantification for hyperbolic and kinetic equations},
  author={Jin, Shi and Pareschi, Lorenzo},
  volume={14},
  year={2018},
  publisher={Springer}
}

@article{HuJin-UQ,
  title={A stochastic {Galerkin} method for the {Boltzmann} equation with uncertainty},
  author={Hu, Jingwei and Jin, Shi},
  journal={Journal of Computational Physics},
  volume={315},
  pages={150--168},
  year={2016},
  publisher={Elsevier}
}

@article{cai2021least,
  title={Least-squares {ReLU neural network (LSNN)} method for linear advection-reaction equation},
  author={Cai, Zhiqiang and Chen, Jingshuang and Liu, Min},
  journal={Journal of Computational Physics},
  pages={110514},
  year={2021}
}

@article{HJJL,
  title={Trend to equilibrium for the kinetic {Fokker-Planck} equation via the neural network approach},
  author={Hwang, Hyung Ju and Jang, Jin Woo and Jo, Hyeontae and Lee, Jae Yong},
  journal={Journal of Computational Physics},
  volume={419},
  pages={109665},
  year={2020}
}

@article{CLM,
  title={Solving the linear transport equation by a deep neural network approach},
  author={Chen, Zheng and Liu, Liu and Mu, Lin},
  journal={arXiv preprint arXiv:2102.09157},
  year={2021}
}

@article{LY,
  title={Asymptotic preserving scheme for anisotropic elliptic equations with deep neural network},
  author={Li, Long and Yang, Chang},
  journal={arXiv preprint arXiv:2104.05337},
  year={2021}
}

@article{lou2021physics,
  title={Physics-informed neural networks for solving forward and inverse flow problems via the {Boltzmann-BGK} formulation},
  author={Lou, Qin and Meng, Xuhui and Karniadakis, George Em},
  journal={Journal of Computational Physics},
  volume={447},
  pages={110676},
  year={2021}
}

@article{DHY,
  title={Multiscale and Nonlocal Learning for {PDEs} using Densely Connected {RNNs}},
  author={Delgadillo, Ricardo A and Hu, Jingwei and Yang, Haizhao},
  journal={arXiv preprint arXiv:2109.01790},
  year={2021}
}

@article{lyu2020mim,
title={{MIM}: A deep mixed residual method for solving high-order partial differential equations},
author={Lyu, Liyao and Zhang, Zhen and Chen, Minxin and Chen, Jingrun},
journal={arXiv preprint arXiv:2006.04146},
year={2020}
}

@article{zang2020weak,
title={Weak adversarial networks for high-dimensional partial differential equations},
author={Zang, Yaohua and Bao, Gang and Ye, Xiaojing and Zhou, Haomin},
journal={Journal of Computational Physics},
pages={109409},
year={2020},
publisher={Elsevier}
}

@Article{deepGalerkin2018,
author={Justin Sirignano and Konstantinos Spiliopoulos},
title={{DGM: A} deep learning algorithm for solving partial differential equations},
journal={Journal of Computational Physics},
volume={375},
year={2018},
pages={1339--1364}
}

@article{beck2020overview,
  title={An overview on deep learning-based approximation methods for partial differential equations},
  author={Beck, Christian and Hutzenthaler, Martin and Jentzen, Arnulf and Kuckuck, Benno},
  journal={arXiv preprint arXiv:2012.12348},
  year={2020}
}

@article{liao2019deep,
title={Deep {N}itsche Method: Deep {R}itz Method with Essential Boundary Conditions},
author={Liao, Yulei and Ming, Pingbing},
journal={arXiv preprint arXiv:1912.01309},
year={2019}
}

@article{li2020fourier,
  title={Fourier neural operator for parametric partial differential equations},
  author={Li, Zongyi and Kovachki, Nikola and Azizzadenesheli, Kamyar and Liu, Burigede and Bhattacharya, Kaushik and Stuart, Andrew and Anandkumar, Anima},
  journal={arXiv preprint arXiv:2010.08895},
  year={2020}
}

@article{lu2021learning,
  title={Learning nonlinear operators via DeepONet based on the universal approximation theorem of operators},
  author={Lu, Lu and Jin, Pengzhan and Pang, Guofei and Zhang, Zhongqiang and Karniadakis, George Em},
  journal={Nature Machine Intelligence},
  volume={3},
  pages={218--229},
  year={2021}
}

@article{raissi2019physics,
title={Physics-informed neural networks: A deep learning framework for solving forward and inverse problems involving nonlinear partial differential equations},
author={Raissi, Maziar and Perdikaris, Paris and Karniadakis, George E},
journal={Journal of Computational Physics},
volume={378},
pages={686--707},
year={2019}
}

@Article{E2018,
author={E, Weinan and Yu, Bing},
title={The Deep {Ritz} Method: A Deep Learning-Based Numerical Algorithm for Solving Variational Problems},
journal={Communications in Mathematics and Statistics},
year={2018},
volume={6},
number={1},
pages={1--12}
}

@article{alexandre2011global,
  title={Global existence and full regularity of the {Boltzmann} equation without angular cutoff},
  author={Alexandre, Radjesvarane and Morimoto, Yoshinori and Ukai, Seiji and Xu, C -J and Yang, Tong},
  journal={Communications in Mathematical Physics},
  volume={304},
  pages={513--581},
  year={2011},
  publisher={Springer}
}

@misc{ukai2007,
  title={Mathematical Theory of {Boltzmann} Equation},
  author={Ukai, Seiji and Yang, Tong},
  howpublished={Lecture Notes, Department of Mathematics and Liu Bie Ju Center for Mathematical Sciences City University of Hong Kong},
  year={2007},
  url={http://www.cityu.edu.hk/rcms/publications/ln8.pdf}
}

\end{document}